\documentclass[reqno]{article}
\usepackage{microtype}
\usepackage{amsmath,amsfonts,amsthm,amssymb}
\usepackage{bm}
\usepackage{indentfirst}
\usepackage{graphicx}
\numberwithin{equation}{section}
\usepackage[margin=2.9cm]{geometry}
\usepackage[colorlinks,linkcolor=blue,citecolor=blue]{hyperref}
\usepackage{xcolor}
\usepackage{amsmath}
\allowdisplaybreaks[4]
\theoremstyle{plain}
\newtheorem{Thm}{Theorem}[section]
\newtheorem*{Thm*}{Theorem}
\newtheorem{Lem}[Thm]{Lemma}

\newtheorem{Prop}[Thm]{Proposition}

\theoremstyle{definition}

\newtheorem{Rem}[Thm]{Remark}

\newcommand{\p}{\partial}

\newcommand{\R}{\mathbb{R}}

\begin{document}

\title{Time-Decay Rates toward Rarefaction Waves for the One-Dimensional Isentropic Compressible Navier--Stokes Equations}

\author{Mutong He\thanks{Academy of Mathematics and Systems Science, Chinese Academy of Sciences, Beijing 100190, PR China. E-mail: hemutong25@mails.ucas.ac.cn},
Feimin Huang\thanks{Academy of Mathematics and Systems Science, Chinese Academy of Sciences, Beijing 100190, PR China. E-mail: fhuang@amt.ac.cn},
Lingjun Liu\thanks{School of Mathematics, Statistics and Mechanics, Beijing University of Technology, Beijing 100124, PR China. E-mail: lingjunliu@bjut.edu.cn},
Lingda Xu\thanks{Hetao Institute of Mathematics and Interdisciplinary Sciences, Shenzhen 518048, PR China. E-mail: xuldmath@gmail.com}
}

\date{}
\maketitle

\begin{abstract}
While the asymptotic stability of rarefaction waves for the compressible Navier--Stokes equations has been well established since the pioneering works of Matsumura--Nishihara \cite{MatsumuraNishihara1986} and Liu--Xin \cite{LiuXin1988}, the classical $L^2$ energy method yields only qualitative convergence in $L^\infty(\R)$ as $t\to\infty$ without any explicit decay rate. In this paper, we establish for the first time a standalone \(L^4\) energy estimate for perturbations of a weak rarefaction wave in the one-dimensional isentropic compressible Navier--Stokes system, thereby obtaining the first explicit algebraic time-decay rates for this model. Specifically, for small initial perturbations in $H^2(\R)$ and small wave strength $\delta$, we prove that $(v-V,u-U)$ decays in $L^4(\R)$ at the rate
\[
\|(v-V,u-U)(t)\|_{L^4}\le C\bigl(\|(v_0-V(0),u_0-U(0))\|_{H^2}+\sqrt{\delta}\bigr)(1+t)^{-1/8},
\]
where \((V,U)\) denotes the smooth rarefaction wave, while the relative Riemann invariants \(z_\pm\) of the perturbed system decay in \(L^\infty(\R)\) at the rate \((1+t)^{-1/6}\). The proof relies on a new local quartic energy estimate for \(z_\pm\) and \(z_\pm-Q\), where \(Q=\mu(v)v_x-\mu(V)V_x\), as well as a nonlocal half-space wave-interaction functional for corrected characteristic modes and a weighted higher-order derivative energy estimate coupled at the scale of the one-dimensional Nash inequality.
\end{abstract}



\tableofcontents

\section{Introduction}\label{1}
We consider the one-dimensional isentropic compressible Navier--Stokes system in Lagrangian coordinates:
\begin{align}\label{eqs1.1}
\begin{cases}
v_t-u_x=0,\\
u_t+P(v)_x=(\mu(v)u_x)_x, \qquad (t,x)\in(0,\infty)\times\R,
\end{cases}
\end{align}
\noindent with initial data
\begin{equation}\label{eqs2.8}
(v,u)|_{t=0}=(v_0,u_0), 
\end{equation}
which satisfies  $(v_0,u_0)(x)\to(v_\pm,u_\pm)$ as $x\to\pm\infty$.
\noindent Here $v>0$ is the specific volume, $u$ is the fluid velocity, $P(v)$ is the pressure satisfying $P'(v)<0$, and $\mu(v)>0$ is the viscosity coefficient. Let $(V,U)(t,x)$ be a smooth profile approximating a weak $2$-rarefaction wave of strength $\delta>0$, and denote the perturbation by
\[
\phi=v-V,\qquad \psi=u-U.
\]
The global existence and large-time behavior of solutions to \eqref{eqs1.1} and related viscous gas dynamics models have been extensively studied since the foundational works of Kanel' \cite{Kanel1968}, Kazhikhov--Shelukhin \cite{KazhikhovShelukhin1977}, and Matsumura--Nishida \cite{MatsumuraNishida1980}. For the Riemann problem with rarefaction wave data, Matsumura--Nishihara \cite{MatsumuraNishihara1986,MatsumuraNishihara1992} and Liu--Xin \cite{LiuXin1988} proved that for small initial perturbations in $H^2(\R)$ and small wave strength $\delta$, \eqref{eqs1.1} admits a unique global strong solution $(v,u)$ converging to $(V,U)$ in $L^\infty(\R)$ as $t\to\infty$; see also \cite{KawashimaMatsumuraNishihara1986,Xin1988,Xin1989,NishiharaYangZhao2004,HuangLiWang2012} for extensions to non-isentropic flows, general $2\times 2$ viscous conservation laws, strong rarefaction waves, and vanishing viscosity limits.

Despite these classical stability results, obtaining an explicit time-decay rate toward rarefaction waves for the compressible Navier--Stokes system \eqref{eqs1.1} has remained a longstanding challenge. Unlike viscous shock waves \cite{MatsumuraNishihara1985} or contact discontinuities \cite{HuangMatsumuraXin2006,HuangXinYang2008}, where anti-derivative variables can be introduced to extract algebraic decay rates, rarefaction waves do not admit an $L^2$ anti-derivative framework. As a result, the classical $H^2$ energy method for rarefaction waves controls only the uniform $H^2(\R)$ norm and the spacetime dissipation $\int_0^\infty \|(\phi_x,\psi_x)(t)\|_{H^1}^2\,dt$, yielding qualitative convergence $\|(\phi,\psi)(t)\|_{L^\infty}\to 0$ as $t\to\infty$ without any quantitative decay rate. For scalar viscous conservation laws, standalone $L^p$ ($p>2$) energy estimates can be derived directly to obtain decay rates toward rarefaction waves via Nash-type inequalities \cite{HattoriNishihara1991,Ito1996}. For the compressible Navier--Stokes system \eqref{eqs1.1}, however, no standalone $L^p$ ($p>2$) estimate was previously available due to the coupled hyperbolic--parabolic structure and the nonlinear coupling between characteristic families.

The main contribution of this paper is to establish, for the first time, a standalone $L^p$ ($p=4$) energy estimate for perturbations of rarefaction waves in the compressible Navier--Stokes system \eqref{eqs1.1}, and consequently to obtain the first explicit algebraic time-decay rates toward rarefaction waves for the compressible Navier--Stokes equations. Specifically, for small initial perturbations in $H^2(\R)$ and small wave strength $\delta$, we establish the $L^4$ decay rate
\begin{align*}
\|(\phi,\psi)(t)\|_{L^4}\le CA(1+t)^{-1/8},\qquad A^2:=\|(\phi,\psi)(0)\|_{H^2}^2+\delta,
\end{align*}
together with the pointwise decay $\|z_\pm(t)\|_{L^\infty}\le CA(1+t)^{-1/6}$ for the relative Riemann invariants $z_\pm$ defined in \eqref{eqs3.1}. Note that $(1+t)^{-1/8}$ coincides with the $L^4$ decay rate of the one-dimensional heat equation with $L^2$ initial data.

Deriving a standalone $L^4$ estimate for the coupled hyperbolic--parabolic system \eqref{eqs1.1} around a rarefaction wave involves two essential difficulties when working in characteristic coordinates. First, because viscosity acts only on the momentum equation, both Riemann invariants $z_\pm$ share the same viscous flux $[\mu(v)u_x-\mu(V)U_x]_x$ with $(u-U)_x\simeq (z_-+z_+)_x$, so multiplying the $z_\pm$ equations by $z_\pm^3$ mixes $(z_-)_x$ and $(z_+)_x$ rather than controlling $z_\pm^2(z_\pm)_x^2$ individually. Second, because the sound speed $w(v)=\sqrt{-P'(v)}$ depends on $v$, integrating the quartic flux terms by parts produces the sixth-order cross integral $\int_\R(z_-^4z_+^2+z_-^2z_+^4)\,dx$. On a rarefaction background, the expansive slope $W_x>0$ of $W=w(V)$ yields a nonnegative damping term proportional to $\int_\R W_x(z_--z_+)^2[(z_-+z_+)^2+2z_+^2]\,dx$ (Proposition~\ref{prop4.1}), which degenerates along the line $z_-=z_+$ and therefore cannot control $\int_\R(z_-^4z_+^2+z_-^2z_+^4)\,dx$.

To close the $L^4$ estimate despite the partial viscosity in characteristic coordinates, we introduce the relative density-gradient variable $Q=\mu(v)v_x-\mu(V)V_x$, motivated by the effective velocity structure in compressible fluid dynamics \cite{Kanel1968,BreschDesjardins2004,Haspot2011}. By the mass equation, $Q_t$ reproduces the viscous term $[\mu(v)u_x-\mu(V)U_x]_x$ up to profile residuals, so subtracting $Q_t$ from the $z_\pm$ equations replaces the second-order viscosity by $\pm w(v)Q_x$. Testing the equations for $z_\pm$ and $z_\pm-Q$ against their cubes and summing the four relations completes two squares and yields coercive control over $\int_\R [z_-^2(z_-)_x^2+z_+^2(z_+)_x^2+Q^4]\,dx$. Moreover, integrating the nonlinear flux in state space before integrating by parts in $x$ cancels the frozen-coefficient quartic profile error, leaving only the mixed sixth-order integral $\int_\R(z_-^4z_+^2+z_-^2z_+^4)\,dx$.

To bound $\int_\R(z_-^4z_+^2+z_-^2z_+^4)\,dx$, we exploit the strict separation of the characteristic speeds $\pm w(v)$ by constructing a nonlocal half-space wave-interaction functional $\mathcal{I}(t)$ inspired by Glimm's interaction potential \cite{Glimm1965}, pairing quartic and quadratic characteristic modes from opposite families. Because $w(v)$ is state-dependent, we modify the pure powers $z_\pm^m$ by a first-order transverse factor $\eta_\pm^{(m)}=z_\pm^m(1\mp\alpha(V)z_\mp)$ so that the compatibility error at $m=2$ becomes higher-order. The resulting cross-family gradient terms $A^2\int_\R z_\mp^2(z_\pm)_x^2\,dx$ generated by viscosity in $\mathcal{I}'(t)$ are controlled after integration by parts by coupling $\mathcal{I}(t)$ with a $W$-weighted higher-order energy $A^2E_1(t)$. By the one-dimensional Nash inequality \cite{Nash1958}, the combined Lyapunov functional $E_{\mathrm{tot}}(t)=E(t)+MA^2E_1(t)+\Lambda\mathcal{I}(t)$ satisfies a closed cubic differential inequality, from which the decay rates follow by ODE comparison.

Let us briefly review related works on time-decay rates toward basic wave patterns. For scalar viscous conservation laws, asymptotic stability and decay rates toward rarefaction and traveling waves were studied in \cite{IlinOleinik1960,HattoriNishihara1991,Ito1996,NishikawaNishihara2000,HuangXu2022}, where scalar $L^p$ estimates and maximum principles are available. For one-dimensional hyperbolic--parabolic systems, pointwise wave structures under $L^1$ or algebraic spatial weights on the initial perturbation were investigated via Green's function methods in \cite{Zeng1994,SzepessyZumbrun1996,LiuZeng1997,LiuYu2022}. For viscous shock waves, contact discontinuities, and composite waves of the compressible Navier--Stokes equations, stability and decay rates via anti-derivative and weighted energy methods were established in \cite{MatsumuraNishihara1985,HuangMatsumuraXin2006,HuangXinYang2008,HuangLiMatsumura2010,HuangXuYuan2022,LiuWangXu2023}. Around constant equilibrium states, decay rates in critical spaces were obtained in \cite{DanchinXu2017,XinXu2021}, and a physical-space $L^p$ approach for partially dissipative systems near constant states was recently developed by Crin-Barat, Shou, and Zuazua \cite{Crin-BaratShouZuazua2025}. In contrast to these works, the present paper establishes the first standalone $L^4$ estimate on a moving rarefaction background and yields the first explicit time-decay rate toward rarefaction waves for the compressible Navier--Stokes system in the $H^2$ framework.

The rest of the paper is organized as follows. Section~\ref{sec2} introduces the smooth rarefaction profile and states the main theorem (Theorem~\ref{thm2.2}). Section~\ref{sec3} derives the characteristic reformulation and establishes the local quartic estimate (Proposition~\ref{prop4.1}). Section~\ref{sec5} proves the weighted higher-order energy estimate and its Nash scaling. Sections~\ref{sec6} and \ref{sec7} construct the corrected characteristic modes and estimate the nonlocal functional $\mathcal{I}(t)$. Section~\ref{sec8} closes the Lyapunov inequality and completes the proof of Theorem~\ref{thm2.2}.

\section{Rarefaction profile and main results}\label{sec2}
Fix a compact interval $I \Subset (0,\infty)$ whose interior contains the range of the background specific volume with a positive margin. Throughout the paper, the perturbation smallness threshold is chosen so that $v(t,x) = V(t,x) + \phi(t,x)\in I$ for all $(t,x)\in[0,\infty)\times\R$. We assume that the pressure $P$, the viscosity $\mu$, and the sound speed $w(v)$ satisfy
\begin{align}\label{eqs2.1}
P\in C^5(I), \qquad \mu\in C^4(I),\qquad w(v):=\sqrt{-P'(v)},
\end{align}
and that there exist constants $c_0, \mu_0, c_g>0$ such that
\begin{align}\label{eqs2.2}
w(v)\ge c_0, \qquad \mu(v)\ge\mu_0,\qquad w'(v)\le-c_g<0 \qquad \forall\,v\in I.
\end{align}
The genuine nonlinearity condition $w'(v)\le-c_g<0$ (equivalently, $P''(v)>0$) fixes the orientation of the $2$-rarefaction wave. Here and below, $W(t,x)=w(V(t,x))$ denotes the background sound speed, and $C,c>0$ denote generic constants independent of $t$ and $\delta$.

\subsection{Smooth rarefaction profile}
Let $(v_\pm,u_\pm)\in\operatorname{int}I\times\R$ be fixed asymptotic end states, and set $W_\pm:=w(v_\pm)$. For a $2$-rarefaction wave connecting $(v_-,u_-)$ to $(v_+,u_+)$, we assume
\begin{align}\label{eqs2.3}
W_-<W_+,\qquad u_+=u_--\int_{v_-}^{v_+}w(y)\,dy,
\end{align}
with $w^{-1}([W_-,W_+])\Subset\operatorname{int}I$. By \eqref{eqs2.2}, $w^{-1}$ is smooth and strictly decreasing on $[W_-,W_+]$.

Following \cite{MatsumuraNishihara1986,HattoriNishihara1991}, let $\widetilde{W}(\tau,x)$ be the classical solution to the viscous Burgers equation
\begin{align*}
\widetilde{W}_\tau+\widetilde{W}\widetilde{W}_x=\widetilde{W}_{xx},\qquad \widetilde{W}(0,x)=\begin{cases}
W_-, & x<0,\\
W_+, & x>0,
\end{cases}
\end{align*}
and define the time-shifted Burgers profile
\begin{align}\label{eqs2.4}
W(t,x):=\widetilde{W}(1+t,x).
\end{align}
Then $W(t,x)$ is smooth for all $t\ge 0$, satisfies $W_t+WW_x=W_{xx}$ with $W_x(t,x)>0$, and converges to $W_\pm$ as $x\to\pm\infty$. We define the smooth rarefaction profile $(V,U)(t,x)$ and its wave strength $\delta$ by
\begin{equation}\label{eqs2.5}
V(t,x):=w^{-1}(W(t,x)),
\qquad
U(t,x):=u_- - \int_{v_-}^{V(t,x)}w(y)\,dy,
\qquad
\delta:=W_+-W_-.
\end{equation}
Differentiating \eqref{eqs2.5} yields
\begin{equation}\label{eqs2.6}
W_x\ge0,
\qquad
U_x=-WV_x,
\qquad
V_x=\frac{W_x}{w'(V)},
\qquad
W_t+WW_x=W_{xx}.
\end{equation}
Substituting $(V,U)$ into \eqref{eqs1.1} produces the profile residuals
\begin{equation}\label{eqs2.7}
R_1:=V_t-U_x=\frac{W_{xx}}{w'(V)},
\qquad
R_2:=U_t+P(V)_x-[\mu(V)U_x]_x.
\end{equation}
The following lemma collects standard \(L^r\) decay estimates for \((V,U)\) and \((R_1,R_2)\); see \cite{kawashima2004stability, LiuXin1988, HattoriNishihara1991, MatsumuraNishihara1992} for details.

\begin{Lem}\label{lem2.1}
Let $0<\delta\le1$. For each $j=0,1,2$ and $1\le r\le\infty$,
\begin{equation}\label{eqs2.9}
\|\p_x^{j+2}W(t)\|_{L^r}
+\|\p_x^j[(W_x)^2](t)\|_{L^r}
+\|\p_x^j(R_1,R_2)(t)\|_{L^r}\le
C_{j,r}(1+t)^{-3/2- j/2+1/(2r)}
\min\{1,\delta\sqrt{1+t}\}.
\end{equation}
Moreover,
\begin{equation}\label{eqs2.10}
W_-<W(t,x)<W_+, \qquad W_x(t,x)\ge0, \qquad \int_{\R}W_x(t,x)\,dx=\delta, \qquad \|W_x(t)\|_{L^\infty} \le C\min\{\delta,(1+t)^{-1}\},
\end{equation}
\begin{equation}\label{eqs2.11}
\Bigl\|\p_x\Bigl(\frac{\mu(V)}{w'(V)}W_x\Bigr)(t)\Bigr\|_{L^r}
\le
C_r(1+t)^{-3/2+1/(2r)}
\min\{1,\delta\sqrt{1+t}\},
\end{equation}
and
\begin{equation}\label{eqs2.12}
\|[\mu(V)R_1]_x(t)\|_{L^r}
\le
C_r(1+t)^{-2+1/(2r)}
\min\{1,\delta\sqrt{1+t}\}.
\end{equation}
In particular,
\begin{equation}\label{eqs2.13}
\int_0^\infty
\Bigl[\|W_{xx}\|_{H^2}
+\|(W_x)^2\|_{H^2}
+\|(R_1,R_2)\|_{H^2}+\Bigl\|\p_x\Bigl(\frac{\mu(V)}{w'(V)}W_x\Bigr)\Bigr\|_{L^2}^2
+\|W_x\|_{L^\infty}^2
+\|W_x\|_{L^\infty}^4
\Bigr]\,dt
\le C\sqrt{\delta},
\end{equation}
\begin{equation}\label{eqs2.14}
\int_0^\infty\|W_x(t)\|_{L^\infty}^{4/3}\,dt
\le C\delta^{1/3},
\qquad
\int_0^\infty
\bigl(\|R_1(t)\|_{L^2}+\|R_2(t)\|_{L^2}\bigr)\,dt
\le C\sqrt{\delta},
\end{equation}
and
\begin{equation}\label{eqs2.15}
\|R_1(t)\|_{L^4}+\|R_2(t)\|_{L^4}
+\|[\mu(V)R_1]_x(t)\|_{L^4}\le
C(1+t)^{-11/8}\min\{1,\delta\sqrt{1+t}\}.
\end{equation}
\end{Lem}

\begin{Rem}
The time-integral bounds in \eqref{eqs2.13}--\eqref{eqs2.14} follow by splitting the integration at $\tau=1+t=\delta^{-2}$ or $\tau=\delta^{-1}$; for example,
\[
\int_1^\infty
\tau^{-5/4}\min\{1,\delta\sqrt{\tau}\}\,d\tau=
\delta\int_1^{\delta^{-2}}\tau^{-3/4}\,d\tau
+\int_{\delta^{-2}}^\infty\tau^{-5/4}\,d\tau\le C\sqrt{\delta}.
\]
\end{Rem}

\subsection{Main results}
For the perturbation $\phi=v-V$ and $\psi=u-U$, we also set
\begin{align}\label{eqs2.16}
h:=\frac{\psi}{W}.
\end{align}

We first state the global existence and \(H^2\) asymptotic stability result used in this paper. This formulation is based on the classical stability theory for rarefaction waves \cite{MatsumuraNishihara1986, LiuXin1988, MatsumuraNishihara1992, NishiharaYangZhao2004}, supplemented by the weighted velocity estimate
\[
\int_0^\infty\int_{\R} W_x h^2\,dx\,dt<+\infty.
\]
\noindent Specifically, by considering the weighted relative entropy
\[
\mathcal E_0(t):=\int_{\R}\frac{W^2}{W_-^4}
\left\{\frac{\psi^2}{2}
-\int_0^\phi\bigl[P(V+s)-P(V)\bigr]\,ds\right\}dx,
\]
we can obtain the classical energy estimate with the additional dissipation. Since the proof is elementary and standard, we omit the details here.

\begin{Thm}[Global existence and $H^2$  stability]\label{thm2.1}
Assume that \eqref{eqs2.1}--\eqref{eqs2.2} hold, and let $(V,U)$ be the smooth rarefaction profile defined in \eqref{eqs2.5}. Let $\left(\phi_0, \psi_0\right):=\left(v_0-V(0, \cdot), u_0-U(0, \cdot)\right) \in H^2(\mathbb{R})$ and set
\begin{align}\label{eqs2.18}
A^2:=\|(\phi_0,\psi_0)\|_{H^2}^2+\delta.
\end{align}
There exist constants $\varepsilon_0>0$ and $K_S>0$, independent of $(\phi_0,\psi_0)$ and $\delta$, such that whenever $A\le \varepsilon_0$, the Cauchy problem \eqref{eqs1.1}, \eqref{eqs2.8} admits a unique global strong solution $(v,u)=(V+\phi,U+\psi)$ satisfying $v(t,x)\in I$ for all $(t,x)\in[0,\infty)\times\R$,
\begin{align}\label{eqs2.17}
(\phi,\psi)\in C([0,\infty);H^2(\R)),\qquad \phi_x\in L^2(0,\infty;H^1(\R)),\qquad \psi_x\in L^2(0,\infty;H^2(\R)),
\end{align}
and the uniform energy estimate
\begin{align}\label{eqs2.19}
\sup_{t\ge0}\|(\phi,\psi)(t)\|_{H^2}^2+\int_0^\infty\Bigl[\|\psi_x(t)\|_{H^2}^2+\|\phi_x(t)\|_{H^1}^2+\int_{\R}W_x\Bigl(\phi^2+\frac{\psi^2}{W^2}\Bigr)\,dx\Bigr]\,dt\le K_SA^2.
\end{align}
In particular, $\lim\limits_{t\to\infty}\|(\phi,\psi)(t)\|_{L^\infty}=0$, and by the one-dimensional Sobolev embedding $H^1(\R)\hookrightarrow L^\infty(\R)$ together with $W_-\le W\le W_+$, it holds that for $t\ge 0$,
\begin{align}\label{eqs2.20}
\|(\phi,h)(t)\|_{H^2}+\|(\phi_x,h_x)(t)\|_{L^\infty}\le CA.
\end{align}
\end{Thm}

Throughout the paper, generic constants $C$ and $c$ may depend on $K_S$, $I$, and the bounds in \eqref{eqs2.1}--\eqref{eqs2.2}, but never on $t$, $A$, or $\delta$. Our main result establishes explicit algebraic time-decay rates for the global solution obtained in Theorem~\ref{thm2.1}.

\begin{Thm}\label{thm2.2}
Assume that \eqref{eqs2.1}--\eqref{eqs2.2} hold, and let $(V,U)$ be the smooth rarefaction profile defined in \eqref{eqs2.5}. There exist constants $0<A_0\le \varepsilon_0$ and $0<\delta_0\le 1$ such that whenever $A\le A_0$ and $\delta\le\delta_0$, the global solution \((v,u)=(V+\phi,U+\psi)\) of \eqref{eqs1.1} and \eqref{eqs2.8} satisfies for \(t\ge 0\) that
\begin{align}\label{eqs2.21}
\|(\phi,\psi)(t)\|_{L^4}\le CA(1+t)^{-1/8}.
\end{align}
Moreover, the relative Riemann invariants \(z_\pm\) defined in \eqref{eqs3.1} satisfy
\begin{align}\label{eqs2.22}
\|z_-(t)\|_{L^\infty}+\|z_+(t)\|_{L^\infty}\le CA(1+t)^{-1/6}.
\end{align}
\end{Thm}

\begin{Rem}
Theorem~\ref{thm2.2} provides the first explicit time-decay rate toward rarefaction waves for the compressible Navier--Stokes system \eqref{eqs1.1} in the $H^2(\R)$ stability framework, obtained by closing a standalone $L^4$ energy estimate. As shown in Section~\ref{sec8}, the proof also yields the derivative decay $\|(\phi_x,\psi_x)(t)\|_{L^2}\le CA(1+t)^{-1/4}$. Moreover, since \eqref{eqs1.1} is invariant under the reflection $\hat v(t,x)=v(t,-x)$, $\hat u(t,x)=-u(t,-x)$, which interchanges the $1$- and $2$-characteristic families, Theorem~\ref{thm2.2} applies equally to weak $1$-rarefaction waves.
\end{Rem}

\section{Local quartic estimates}\label{sec3}
\subsection{Characteristic reformulation}
Subtracting \eqref{eqs2.7} from \eqref{eqs1.1} gives the perturbation system
\begin{align*}
\begin{cases}
\phi_t-\psi_x=-R_1,\\
\psi_t-w(v)^2\phi_x=(\mu(v)u_x-\mu(V)U_x)_x-R_2-(P'(v)-P'(V))V_x.
\end{cases}
\end{align*}
We define the relative Riemann invariants $z_\pm$ and the relative density-gradient variable $Q$ by
\begin{align}\label{eqs3.1}
z_-:=\psi+\int_V^vw(y)\,dy,\qquad z_+:=\psi-\int_V^vw(y)\,dy,\qquad Q:=\mu(v)v_x-\mu(V)V_x.
\end{align}
By $v_t=u_x$ and $V_t=U_x+R_1$, we have $(\mu(v)v_x)_t=[\mu(v)u_x]_x$ and $(\mu(V)V_x)_t=[\mu(V)(U_x+R_1)]_x$, so that
\begin{align}\label{eqs3.2}
Q_t=[\mu(v)u_x-\mu(V)U_x]_x-[\mu(V)R_1]_x.
\end{align}
Subtracting \eqref{eqs3.2} from the equations for $z_\pm$ eliminates the second-order viscous term from $z_\pm-Q$ and replaces it with $\pm w(v)Q_x$:
\begin{align}
&(z_-)_t-w(v)(z_-)_x=[\mu(v)u_x-\mu(V)U_x]_x-(R_2+WR_1), \label{eqs3.3}\\
&(z_+)_t+w(v)(z_+)_x=[\mu(v)u_x-\mu(V)U_x]_x-(R_2-WR_1)+2W(w(v)-W)V_x,\label{eqs3.4}\\
&(z_--Q)_t-w(v)(z_--Q)_x=w(v)Q_x-\bigl(R_2-[\mu(V)R_1]_x+WR_1\bigr),  \label{eqs3.5}\\
&(z_+-Q)_t+w(v)(z_+-Q)_x=-w(v)Q_x-\bigl(R_2-[\mu(V)R_1]_x-WR_1\bigr)+2W(w(v)-W)V_x. \label{eqs3.6}
\end{align}

We define the local quartic energy
\begin{align}\label{eqs3.7}
E(t):=\frac14\int_{\R}\bigl(z_-^4+z_+^4+(z_--Q)^4+(z_+-Q)^4\bigr)\,dx.
\end{align}
Since
\begin{align}\label{exact}
z_-+z_+=2\psi,\qquad z_--z_+=2\int_{V}^{V+\phi}w(y)\,dy,
\end{align}
and $|a|^4+|a-b|^4\simeq|a|^4+|b|^4$, the bounds $c_0\le w(v)\le C$ give
\begin{align}\label{eqs3.8}
c\int_{\R}\bigl(\phi^4+h^4+Q^4\bigr)\,dx\le E(t)\le C\int_{\R}\bigl(\phi^4+h^4+Q^4\bigr)\,dx.
\end{align}

\subsection{Quartic energy estimate}\label{sec4}

\begin{Prop}\label{prop4.1}
For sufficiently small $A$ and $\delta$, the quartic energy $E(t)$ satisfies
\begin{align}\label{eqs4.1}
&E'(t)+c\int_{\R}\bigl[z_-^2(z_-)_x^2+z_+^2(z_+)_x^2+Q^4\bigr]\,dx+c\int_{\R} W_x(z_--z_+)^2\bigl[(z_-+z_+)^2+2z_+^2\bigr]\,dx\nonumber\\
&\le C\int_{\R}(z_-^4z_+^2+z_-^2z_+^4)\,dx+C\bigl[\|W_x\|_{L^\infty}(\|z_-\|_{L^\infty}+\|z_+\|_{L^\infty})+\|W_x\|^2_{L^\infty}+\|W_x\|_{L^\infty}^4\bigr]E(t)\nonumber\\
&\quad+CE(t)^{3/4}\bigl(\|R_1\|_{L^4}+\|R_2\|_{L^4}+\|[\mu(V)R_1]_x\|_{L^4}\bigr).
\end{align}
Moreover, the time-dependent coefficients multiplying $E(t)$ on the right-hand side of \eqref{eqs4.1} belong to $L^1(0,\infty)$ by Lemma~\ref{lem2.1} and Theorem~\ref{thm2.1}.
\end{Prop}

\begin{proof}
\noindent\textbf{Step 1.}
Using $U_x=-WV_x$, $Q=\mu(v)v_x-\mu(V)V_x$, and $(z_\pm)_x=u_x-U_x\mp w(v)(v_x-V_x)\mp(w(v)-W)V_x$, we decompose the viscous flux as
\begin{align}\label{eqs4.2}
\mu(v)u_x-\mu(V)U_x+w(v)Q&=\mu(v)(z_-)_x+\mu(V)(W-w(v))V_x,\nonumber\\
\mu(v)u_x-\mu(V)U_x-w(v)Q&=\mu(v)(z_+)_x+\bigl[(W+w(v))\mu(V)-2W\mu(v)\bigr]V_x.
\end{align}
Multiplying \eqref{eqs3.3}--\eqref{eqs3.6} by
\begin{align}\label{eqs3.9}
z_-^3,\qquad z_+^3,\qquad (z_--Q)^3,\qquad (z_+-Q)^3,
\end{align}
respectively, and integrating by parts over $\R$ yields
\begin{align*}
&\frac{1}{4}\frac{d}{dt}\int_{\R} z_-^4\,dx+\frac{1}{4}\int_{\R} w(v)_x z_-^4\,dx
=-3\int_{\R}[\mu(v)u_x-\mu(V)U_x]z_-^2(z_-)_x\,dx-\int_{\R} z_-^3(R_2+WR_1)\,dx,\\
&\frac{1}{4}\frac{d}{dt}\int_{\R} z_+^4\,dx-\frac{1}{4}\int_{\R} w(v)_x z_+^4\,dx
=-3\int_{\R}[\mu(v)u_x-\mu(V)U_x]z_+^2(z_+)_x\,dx-\int_{\R} z_+^3(R_2-WR_1)\,dx\\
&\quad+2\int_{\R} W(w(v)-W)V_x z_+^3\,dx,\\
&\frac{1}{4}\frac{d}{dt}\int_{\R}(z_--Q)^4\,dx
=\int_{\R} w(v)(z_--Q)^3(z_-)_x\,dx-\int_{\R}(z_--Q)^3\bigl(R_2-[\mu(V)R_1]_x+WR_1\bigr)\,dx,\\
&\frac{1}{4}\frac{d}{dt}\int_{\R}(z_+-Q)^4\,dx
=-\int_{\R} w(v)(z_+-Q)^3(z_+)_x\,dx+2\int_{\R} W(w(v)-W)V_x(z_+-Q)^3\,dx\\
&\quad-\int_{\R}(z_+-Q)^3\bigl(R_2-[\mu(V)R_1]_x-WR_1\bigr)\,dx.
\end{align*}

Summing these four identities, expanding $(z_\pm-Q)^3$, substituting \eqref{eqs4.2}, and writing
\[
w(v)_x=w'(v)V_x+\frac{w'(v)}{2w(v)}\bigl[(z_-)_x-(z_+)_x-2(w(v)-W)V_x\bigr],
\]
we complete the squares in $z_-(z_-)_x$ and $z_+(z_+)_x$ to obtain
\begin{align}\label{eqs4.3}
&E'(t)+\int_{\R}\Biggl[3\mu(v)\Bigl(z_-(z_-)_x-\frac{w(v)Q^2}{2\mu(v)}\Bigr)^2+3\mu(v)\Bigl(z_+(z_+)_x+\frac{w(v)Q^2}{2\mu(v)}\Bigr)^2+\frac{w(v)^2}{2\mu(v)}Q^4\Biggr]\,dx\nonumber\\
&\quad+\int_{\R} V_x\Bigl[\frac{Ww'(v)}{2w(v)}(z_-^4-z_+^4)-4W(w(v)-W)z_+^3\Bigr]\,dx\nonumber\\
&=\frac14\int_{\R}\frac{w'(v)}{w(v)}\bigl((z_-)_x-(z_+)_x\bigr)(z_+^4-z_-^4)\,dx-3\int_{\R}\mu(V)(W-w(v))V_xz_-^2(z_-)_x\,dx\nonumber\\
&\quad-3\int_{\R}\bigl[(W+w(v))\mu(V)-2W\mu(v)\bigr]V_xz_+^2(z_+)_x\,dx+2\int_{\R}\frac{w(v)}{\mu(v)}\bigl[W\mu(v)-w(v)\mu(V)\bigr]V_xQ^3\,dx\nonumber\\
&\quad-\int_{\R} 2W(w(v)-W)V_x(3z_+^2Q-3z_+Q^2+Q^3)\,dx
-\int_{\R}\Bigl[z_-^3(R_2+WR_1)+z_+^3(R_2-WR_1)\nonumber\\
&\quad+(z_--Q)^3\bigl(R_2-[\mu(V)R_1]_x+WR_1\bigr)+(z_+-Q)^3\bigl(R_2-[\mu(V)R_1]_x-WR_1\bigr)\Bigr]\,dx
=:\sum_{i=1}^6 H_i.
\end{align}
By Young's inequality and the uniform bounds on $w(v)$ and $\mu(v)$,
\begin{align}\label{eqs3.12}
&\int_{\R}\Biggl[3\mu(v)\Bigl(z_-(z_-)_x-\frac{w(v)Q^2}{2\mu(v)}\Bigr)^2+3\mu(v)\Bigl(z_+(z_+)_x+\frac{w(v)Q^2}{2\mu(v)}\Bigr)^2+\frac{w(v)^2}{2\mu(v)}Q^4\Biggr]\,dx\nonumber\\
&\ge c\int_{\R}\bigl[z_-^2(z_-)_x^2+z_+^2(z_+)_x^2+Q^4\bigr]\,dx.
\end{align}

\medskip
\noindent\textbf{Step 2.}
Since $z_--z_+=2\int_V^v w(y)\,dy=2W(v-V)+O(|v-V|^2)$, Taylor's theorem gives
\begin{align*}
\Bigl|\frac{Ww'(v)}{w(v)w'(V)}-1\Bigr|+\Bigl|\frac{2W(w(v)-W)}{w'(V)(z_--z_+)}-1\Bigr|\le C|z_--z_+|,
\end{align*}
where the second quotient takes its continuous limit $1$ along $z_-=z_+$. Using $V_x=W_x/w'(V)$ and
\[
\frac12(z_-^4-z_+^4)-2(z_--z_+)z_+^3=\frac12(z_--z_+)^2\bigl[(z_-+z_+)^2+2z_+^2\bigr],
\]
we rewrite the integrand of the second term on the left-hand side of \eqref{eqs4.3} as
\begin{align*}
&V_x\Bigl[\frac{Ww'(v)}{2w(v)}(z_-^4-z_+^4)-4W(w(v)-W)z_+^3\Bigr]
=\frac{W_x}{2}(z_--z_+)^2\bigl[(z_-+z_+)^2+2z_+^2\bigr]\\
&\quad+\frac{W_x}{2}\Bigl(\frac{Ww'(v)}{w(v)w'(V)}-1\Bigr)(z_-^4-z_+^4)-2W_x\Bigl(\frac{2W(w(v)-W)}{w'(V)(z_--z_+)}-1\Bigr)(z_--z_+)z_+^3.
\end{align*}
Since both remainder terms contain a factor of $(z_--z_+)$ and their Taylor coefficients are bounded by $C|z_--z_+|\le CA$, we obtain
\begin{equation}\label{eqs4.4}
\int_{\R} V_x\Bigl[\frac{Ww'(v)}{2w(v)}(z_-^4-z_+^4)-4W(w(v)-W)z_+^3\Bigr]\,dx
\ge\frac{1-CA}{2}\int_{\R} W_x(z_--z_+)^2\bigl[(z_-+z_+)^2+2z_+^2\bigr]\,dx.
\end{equation}

\medskip
\noindent\textbf{Step 3.}
Write
\begin{align*}
H_5&=-6\int_{\R} W(w(v)-W)V_x z_+^2 Q\,dx+6\int_{\R} W(w(v)-W)V_x z_+ Q^2\,dx\\
&\quad-2\int_{\R} W(w(v)-W)V_x Q^3\,dx=:\sum_{i=1}^3 H_{5i}.
\end{align*}
Since $|w(v)-W|\le C|z_--z_+|$ and $|V_x|\le CW_x$, applying Young's inequality twice for any $\beta>0$ gives
\begin{align*}
|H_{51}|&\le C\int_{\R} W_x|z_--z_+|z_+^2|Q|\,dx\\
&\le \beta\int_{\R} W_x(z_--z_+)^2\bigl[(z_-+z_+)^2+2z_+^2\bigr]\,dx+C_\beta\int_{\R} W_x z_+^2 Q^2\,dx\\
&\le \beta\int_{\R} W_x(z_--z_+)^2\bigl[(z_-+z_+)^2+2z_+^2\bigr]\,dx+\beta\int_{\R} Q^4\,dx+C_\beta\|W_x\|_{L^\infty}^2 E(t),
\end{align*}
and similarly
\begin{align*}
|H_{52}|+|H_{53}| \le 2\beta\int_{\R} Q^4\,dx+C_\beta\bigl(\|W_x\|_{L^\infty}^2+\|W_x\|_{L^\infty}^4\bigr) E(t).
\end{align*}
Thus
\begin{equation}\label{eqs4.5}
|H_5|\le 3\beta\int_{\R} Q^4\,dx+\beta\int_{\R} W_x(z_--z_+)^2\bigl[(z_-+z_+)^2+2z_+^2\bigr]\,dx+C_\beta\bigl(\|W_x\|_{L^\infty}^2+\|W_x\|_{L^\infty}^4\bigr)E(t).
\end{equation}
For $H_2, H_3, H_4$ in \eqref{eqs4.3}, each coefficient vanishes at $v=V$, so
\begin{align*}
|\mu(V)(W-w(v))|+\bigl|(W+w(v))\mu(V)-2W\mu(v)\bigr|+|W\mu(v)-w(v)\mu(V)|\le C|z_--z_+|.
\end{align*}
By Cauchy--Schwarz and Young's inequalities,
\begin{align*}
|H_2|+|H_3|+|H_4|&\le C\int_{\R} W_x|z_--z_+|\bigl(z_-^2|(z_-)_x|+z_+^2|(z_+)_x|+|Q|^3\bigr)\,dx\\
&\le \beta\int_{\R}\bigl[z_-^2(z_-)_x^2+z_+^2(z_+)_x^2+Q^4\bigr]\,dx+C_\beta\bigl(\|W_x\|_{L^\infty}^2+\|W_x\|_{L^\infty}^4\bigr)E(t).
\end{align*}
By H\"older's inequality and \eqref{eqs3.7}, $H_6$ satisfies
\begin{align}\label{eqs4.7}
|H_6|&\le C\bigl(\|z_-\|_{L^4}^3+\|z_+\|_{L^4}^3+\|z_--Q\|_{L^4}^3+\|z_+-Q\|_{L^4}^3\bigr)\bigl(\|R_1\|_{L^4}+\|R_2\|_{L^4}+\|[\mu(V)R_1]_x\|_{L^4}\bigr)\nonumber\\
&\le CE(t)^{3/4}\bigl(\|R_1\|_{L^4}+\|R_2\|_{L^4}+\|[\mu(V)R_1]_x\|_{L^4}\bigr).
\end{align}

\medskip
\noindent\textbf{Step 4.}
We next show that for every $\beta > 0$,
\begin{align}\label{eqs4.8}
|H_1|&=\Bigl|\frac14\int_{\R}\frac{w'(v)}{w(v)}\bigl((z_-)_x-(z_+)_x\bigr)(z_+^4-z_-^4)\,dx\Bigr|\nonumber\\
&\le\beta\int_{\R}\bigl[z_-^2(z_-)_x^2+z_+^2(z_+)_x^2+Q^4\bigr]\,dx+C_\beta\int_{\R}(z_-^4z_+^2+z_-^2z_+^4)\,dx\nonumber\\
&\quad+C\|W_x(t)\|_{L^\infty}\bigl(\|z_-(t)\|_{L^\infty}+\|z_+(t)\|_{L^\infty}\bigr)E(t).
\end{align}
Define the state parameterization $v(\xi;V)$ implicitly by $2\int_V^{v(\xi;V)}w(y)\,dy=\xi$, so that $v=v(z_--z_+;V)$ and
\begin{align*}
\p_\xi w(v(\xi;V))=\frac12\frac{w'(v(\xi;V))}{w(v(\xi;V))},\qquad \p_Vw(v(\xi;V))=W\frac{w'(v(\xi;V))}{w(v(\xi;V))}.
\end{align*}
By the chain rule,
\begin{align*}
\frac14\frac{w'(v)}{w(v)}\bigl((z_-)_x-(z_+)_x\bigr)=\frac12 w(v)_x-\frac{W}{2}\frac{w'(v)}{w(v)}V_x.
\end{align*}
Integrating $\frac12\int_\R w(v)_x(z_+^4-z_-^4)\,dx$ by parts yields $2\int_\R w(v)[z_-^3(z_-)_x-z_+^3(z_+)_x]\,dx$. Expanding $w(v(z_--z_+;V))$ around $z_+=0$ and $z_-=0$ gives
\begin{align*}
w(v(z_--z_+;V))&=w(v(z_-;V))-\frac12\int_0^{z_+}\frac{w'(v(z_--r;V))}{w(v(z_--r;V))}\,dr,\\
w(v(z_--z_+;V))&=w(v(-z_+;V))+\frac12\int_0^{z_-}\frac{w'(v(r-z_+;V))}{w(v(r-z_+;V))}\,dr.
\end{align*}
Substituting these expansions into $2w(v)[z_-^3(z_-)_x-z_+^3(z_+)_x]$ and integrating the single-mode terms $2w(v(z_-;V))z_-^3(z_-)_x$ and $-2w(v(-z_+;V))z_+^3(z_+)_x$ by parts in $x$, we obtain
\begin{align}\label{eqs4.9}
&\frac14\int_{\R} \frac{w'(v)}{w(v)}\bigl((z_-)_x-(z_+)_x\bigr)(z_+^4-z_-^4)\,dx\nonumber\\
&=-\int_{\R} z_-^3\Bigl[\int_0^{z_+}\frac{w'(v(z_--r;V))}{w(v(z_--r;V))}\,dr\Bigr](z_-)_x\,dx-\int_{\R} z_+^3\Bigl[\int_0^{z_-}\frac{w'(v(r-z_+;V))}{w(v(r-z_+;V))}\,dr\Bigr](z_+)_x\,dx\nonumber\\
&\quad-\int_{\R} V_x\Biggl[2W\int_0^{z_-}\frac{w'(v(r;V))}{w(v(r;V))}r^3\,dr-2W\int_0^{z_+}\frac{w'(v(-r;V))}{w(v(-r;V))}r^3\,dr\nonumber\\
&\qquad\qquad\quad+\frac{W}{2}\frac{w'(v(z_--z_+;V))}{w(v(z_--z_+;V))}(z_+^4-z_-^4)\Biggr]\,dx.
\end{align}
Since $|w'/w|\le C$ on $I$, the state integrals in the first two terms are bounded by $C|z_+|$ and $C|z_-|$, respectively, so Cauchy--Schwarz and Young's inequalities bound the first two integrals in \eqref{eqs4.9} by
\begin{align*}
\beta\int_{\R} \bigl[z_-^2(z_-)_x^2+z_+^2(z_+)_x^2\bigr]\,dx+C_\beta\int_{\R} (z_-^4z_+^2+z_-^2z_+^4)\,dx.
\end{align*}
In the last integral of \eqref{eqs4.9}, freezing $w'/w$ inside both $\int_0^{z_\pm}(\cdots)r^3\,dr$ integrals at $v(z_--z_+;V)$ cancels the quartic terms identically:
\begin{align*}
\frac{W}{2}\frac{w'(v(z_--z_+;V))}{w(v(z_--z_+;V))}\bigl(z_-^4-z_+^4+z_+^4-z_-^4\bigr)=0.
\end{align*}
By the mean-value theorem and $|V_x|\le CW_x$, the last integral in \eqref{eqs4.9} is therefore bounded by
\begin{align*}
C\|W_x(t)\|_{L^\infty}\bigl(\|z_-(t)\|_{L^\infty}+\|z_+(t)\|_{L^\infty}\bigr)E(t),
\end{align*}
which gives \eqref{eqs4.8}.

\medskip
\noindent\textbf{Step 5.}
Taking $A$ small enough that $1-CA\ge1/2$ in \eqref{eqs4.4}, choosing $\beta>0$ sufficiently small, and combining \eqref{eqs4.3}--\eqref{eqs4.8} yields \eqref{eqs4.1}.

To verify that $\|W_x(t)\|_{L^\infty}(\|z_-(t)\|_{L^\infty}+\|z_+(t)\|_{L^\infty})\in L^1(0,\infty)$, we differentiate \eqref{eqs3.1}:
\begin{align}\label{eqs4.10}
(z_-)_x=\psi_x+w(v)\phi_x+(w(v)-W)V_x,\qquad (z_+)_x=\psi_x-w(v)\phi_x-(w(v)-W)V_x.
\end{align}
By \eqref{eqs2.10} and \eqref{eqs2.19},
\[
\int_0^\infty\|(w(v)-W)V_x\|_{L^2}^2\,dt\le C\int_0^\infty\|W_x(t)\|_{L^\infty}\Bigl(\int_\R W_x\phi^2\,dx\Bigr)\,dt\le C\delta A^2,
\]
and hence $\int_0^\infty\|(z_-,z_+)_x(t)\|_{L^2}^2\,dt\le CA^2$. By the Gagliardo--Nirenberg inequality $\|z_\pm(t)\|_{L^\infty}\le CA^{1/2}\|(z_\pm)_x(t)\|_{L^2}^{1/2}$, H\"older's inequality in time with exponents $(4/3,4)$, and \eqref{eqs2.14}, we obtain
\begin{align}\label{eqs4.12}
&\int_0^\infty\|W_x(t)\|_{L^\infty}\bigl(\|z_-(t)\|_{L^\infty}+\|z_+(t)\|_{L^\infty}\bigr)\,dt\nonumber\\
&\le CA^{1/2}\Bigl(\int_0^\infty\|W_x(t)\|_{L^\infty}^{4/3}\,dt\Bigr)^{3/4}\Bigl(\int_0^\infty\|(z_-,z_+)_x(t)\|_{L^2}^2\,dt\Bigr)^{1/4}\le CA\delta^{1/4}\le CA^{3/2}.
\end{align}
\end{proof}

\section{Weighted higher-order energy estimates}\label{sec5}
To control the second-order derivatives $\|(z_\pm)_{xx}\|_{L^2}^2$ needed in Section~\ref{sec7}, we rewrite the perturbation system in terms of $\phi=v-V$ and $h=\psi/W$:
\begin{align}\label{eqs5.1}
\begin{cases}
\displaystyle\phi_t-Wh_x=W_xh-R_1,\\[4pt]
\displaystyle h_t-\frac{w(V+\phi)^2}{W}\phi_x=[\mu(V+\phi)h_x]_x+\p_x(W_x g)+\frac{W_x\mu(V+\phi)}{W}h_x\\[6pt]
\displaystyle\qquad\qquad\qquad\qquad\quad+W_x p+\frac{W_x^2}{W}g-\frac{W_{xx}}{W}h-\frac{R_2}{W},
\end{cases}
\end{align}
where
\begin{align*}
g:=\frac{\mu(V+\phi)}{W}h-\frac{\mu(V+\phi)-\mu(V)}{w'(V)},\qquad p:=h-\frac{P'(V+\phi)-P'(V)}{Ww'(V)}.
\end{align*}
With the profile weight
\begin{align}\label{eqs5.2}
\omega:=\left(\frac{W}{W_-}\right)^4
\end{align}
and a small constant $\kappa>0$, we define the differentiated energy functional $E_1(t)$ and its associated dissipation $D_1(t)$ by
\begin{align}\label{eqs5.3}
E_1(t)&:=\frac12\int_{\R}\omega\frac{w(V+\phi)^2}{W^2}(\phi_x^2+\phi_{xx}^2)\,dx+\frac12\int_{\R}\omega(h_x^2+h_{xx}^2)\,dx-\kappa\int_{\R}\omega h_x\phi_{xx}\,dx,\\
D_1(t)&:=\|\phi_{xx}\|_{L^2}^2+\|h_{xx}\|_{L^2}^2+\|h_{xxx}\|_{L^2}^2+\int_{\R} W_x(\phi_x^2+h_x^2+\phi_{xx}^2+h_{xx}^2)\,dx,\label{def:D1}
\end{align}
and set
\begin{align}\label{eqs5.4}
\mathcal B(t)&:=\|\psi_x\|_{H^2}^2+\|\phi_x\|_{H^1}^2+\int_{\R} W_x(\phi^2+\psi^2)\,dx+\|W_{xx}\|_{H^2}\nonumber\\
&\quad+\|(W_x)^2\|_{H^2}+\|(R_1,R_2)\|_{H^2}+\Bigl\|\p_x\Bigl(\frac{\mu(V)}{w'(V)}W_x\Bigr)\Bigr\|^2_{L^2}.
\end{align}

\begin{Lem}\label{lem5.1}
For $\kappa>0$ fixed sufficiently small, $E_1(t)$ satisfies
\begin{align}\label{eqs5.5}
c\|(\phi_x,h_x)(t)\|_{H^1}^2\le E_1(t)\le C\|(\phi_x,h_x)(t)\|_{H^1}^2,\qquad E_1(0)\le CA^2.
\end{align}
Moreover, for every $\beta>0$ and sufficiently small $A_0,\delta_0$,
\begin{align}\label{eqs5.6}
A^2E'_1(t)+cA^2D_1(t)
&\le \beta\int_{\R}\bigl[z_-^2(z_-)_x^2+z_+^2(z_+)_x^2+Q^4\bigr]\,dx+\beta\int_{\R} W_x(z_--z_+)^2\bigl[(z_-+z_+)^2+2z_+^2\bigr]\,dx\nonumber\\
&\quad+\beta\int_{\R}(z_-^4z_+^2+z_-^2z_+^4)\,dx+C_\beta\mathcal B(t)\bigl(E(t)+A^2E_1(t)\bigr)+CA^4(1+t)^{-3/2}.
\end{align}
\end{Lem}

\begin{proof}
\noindent\textbf{Step 1.}
For $m\in\{1,2\}$, applying $\p_x^m$ to \eqref{eqs5.1}, multiplying the two equations by $\omega\frac{w(V+\phi)^2}{W^2}\p_x^m\phi$ and $\omega\p_x^mh$, respectively, and integrating by parts over $\R$ gives
\begin{align}\label{eqs5.7}
&\frac{d}{dt}\frac12\int_{\R}\omega\Bigl(\frac{w(V+\phi)^2}{W^2}|\p_x^m\phi|^2+|\p_x^mh|^2\Bigr)\,dx+\int_{\R}\omega\mu(V+\phi)|\p_x^{m+1}h|^2\,dx\nonumber\\
&=\frac12\int_{\R}\p_t\Bigl(\omega\frac{w(V+\phi)^2}{W^2}\Bigr)|\p_x^m\phi|^2\,dx+\frac12\int_{\R}\p_t\omega\,|\p_x^mh|^2\,dx-\int_{\R}\p_x\Bigl(\omega\frac{w(V+\phi)^2}{W}\Bigr)\p_x^m\phi\,\p_x^mh\,dx\nonumber\\
&\quad+\int_{\R}\omega\frac{w(V+\phi)^2}{W^2}\p_x^m\phi\bigl[\p_x^m(Wh_x)-W\p_x^{m+1}h+\p_x^m(W_xh-R_1)\bigr]\,dx\nonumber\\
&\quad+\int_{\R}\omega\p_x^mh\Bigl[\p_x^m\Bigl(\frac{w(V+\phi)^2}{W}\phi_x\Bigr)-\frac{w(V+\phi)^2}{W}\p_x^{m+1}\phi\Bigr]\,dx\nonumber\\
&\quad+\int_{\R}\omega\p_x^mh\,\p_x^m\Bigl(\frac{W_x\mu(V+\phi)}{W}h_x+W_x p+\frac{W_x^2}{W}g-\frac{W_{xx}}{W}h-\frac{R_2}{W}\Bigr)\,dx\nonumber\\
&\quad-\int_{\R}(\omega\p_x^mh)_x\bigl[\p_x^m(\mu(V+\phi)h_x)-\mu(V+\phi)\p_x^{m+1}h+\p_x^m(W_x g)\bigr]\,dx\nonumber\\
&\quad-\int_{\R}\p_x\omega\,\mu(V+\phi)\p_x^mh\,\p_x^{m+1}h\,dx.
\end{align}

\medskip
\noindent\textbf{Step 2.}
To recover dissipation for $\phi_{xx}$, we differentiate the momentum equation in \eqref{eqs5.1} once in $x$ and multiply by $-\omega\phi_{xx}$, differentiate the mass equation twice in $x$ and multiply by $-\omega h_x$, and sum the integrated results:
\begin{align}\label{eqs5.8}
&\frac{d}{dt}\Bigl(-\int_{\R}\omega h_x\phi_{xx}\,dx\Bigr)+\int_{\R}\omega\frac{w(V+\phi)^2}{W}\phi_{xx}^2\,dx\nonumber\\
&=\int_{\R}\omega W h_{xx}^2\,dx-\int_{\R}\omega\mu(V+\phi)h_{xxx}\phi_{xx}\,dx-\int_{\R}\p_t\omega\,h_x\phi_{xx}\,dx-\int_{\R}\omega\p_x\Bigl(\frac{w(V+\phi)^2}{W}\Bigr)\phi_x\phi_{xx}\,dx\nonumber\\
&\quad-\int_{\R}\omega\bigl(2[\mu(V+\phi)]_xh_{xx}+[\mu(V+\phi)]_{xx}h_x\bigr)\phi_{xx}\,dx-\int_{\R}\omega\p_x^2(W_x g)\phi_{xx}\,dx\nonumber\\
&\quad-\int_{\R}\omega\p_x\Bigl(\frac{W_x\mu(V+\phi)}{W}h_x+W_x p+\frac{W_x^2}{W}g-\frac{W_{xx}}{W}h-\frac{R_2}{W}\Bigr)\phi_{xx}\,dx\nonumber\\
&\quad-\int_{\R}\omega h_x(W_xh-R_1)_{xx}\,dx-\int_{\R}\omega\Bigl(\frac{6W_x^2}{W}+\frac52 W_{xx}\Bigr)h_x^2\,dx,
\end{align}
where we used $\p_x\omega=4\omega W_x/W$ to integrate $-\int_\R \omega h_x(Wh_x)_{xx}\,dx$ by parts:
\begin{align*}
-\int_{\R}\omega h_x(Wh_x)_{xx}\,dx=\int_{\R}\omega Wh_{xx}^2\,dx-\int_{\R}\omega\Bigl(\frac{6W_x^2}{W}+\frac52 W_{xx}\Bigr)h_x^2\,dx.
\end{align*}
At $m = 2$, the viscous terms in \eqref{eqs5.7} expand as
\begin{align}\label{eqs5.9}
\int_{\R}\omega h_{xx}\p_x^2[\mu(V+\phi)h_x]_x\,dx&=-\int_{\R}\omega\mu(V+\phi)h_{xxx}^2\,dx-2\int_{\R}\omega[\mu(V+\phi)]_xh_{xx}h_{xxx}\,dx\nonumber\\
&\quad-\int_{\R}\p_x\omega\,\mu(V+\phi)h_{xx}h_{xxx}\,dx-\int_{\R}\omega[\mu(V+\phi)]_{xx}h_xh_{xxx}\,dx\nonumber\\
&\quad-2\int_{\R}\p_x\omega\,[\mu(V+\phi)]_xh_{xx}^2\,dx-\int_{\R}\p_x\omega\,[\mu(V+\phi)]_{xx}h_xh_{xx}\,dx,
\end{align}
and the $m=2$ hyperbolic commutator in \eqref{eqs5.7} contains no $\phi_{xxx}$ term since
\begin{align*}
\p_x^2\Bigl(\frac{w(V+\phi)^2}{W}\phi_x\Bigr)-\frac{w(V+\phi)^2}{W}\phi_{xxx}=2\p_x\Bigl(\frac{w(V+\phi)^2}{W}\Bigr)\phi_{xx}+\p_x^2\Bigl(\frac{w(V+\phi)^2}{W}\Bigr)\phi_x.
\end{align*}
Summing \eqref{eqs5.7} over $m=1,2$ and adding $\kappa$ times \eqref{eqs5.8} gives the evolution of $E_1(t)$.

\medskip
\noindent\textbf{Step 3.}
From $W_t+WW_x=W_{xx}$, we have
\begin{align*}
\p_t\omega=4\omega\Bigl(\frac{W_{xx}}{W}-W_x\Bigr), \qquad \p_x\omega=\frac{4\omega W_x}{W}, \qquad \p_x(\omega W)=5\omega W_x.
\end{align*}
Since $P''(V)=-2Ww'(V)$, Taylor expansion gives $W_x p = W_x(h+2\phi)+O(W_x\phi^2)$. Freezing the state-dependent coefficients in \eqref{eqs5.7} at $\phi=h=0$ (modulo $W_{xx}$-terms), we see that:
\begin{itemize}
\item $\frac12\p_t\omega(|\p_x^m\phi|^2+|\p_x^mh|^2)$ contributes $-2\omega W_x(|\p_x^m\phi|^2+|\p_x^mh|^2)$;
\item $-\p_x(\omega W)\p_x^m\phi\,\p_x^mh$ contributes $-5\omega W_x\p_x^m\phi\,\p_x^mh$;
\item $\p_x^m(Wh_x)-W\p_x^{m+1}h+\p_x^m(W_xh)$ contributes $(m+1)\omega W_x\p_x^m\phi\,\p_x^mh$;
\item $\p_x^m(W\phi_x)-W\p_x^{m+1}\phi$ and $\p_x^m[W_x(h+2\phi)]$ contribute $(m+2)\omega W_x\p_x^m\phi\,\p_x^mh+\omega W_x|\p_x^mh|^2$.
\end{itemize}
Moving the sum of these terms to the left-hand side of \eqref{eqs5.7} yields
\begin{align*}
\int_{\R}\omega W_x\bigl[2|\p_x^m\phi|^2+|\p_x^mh|^2-(2m-2)\p_x^m\phi\,\p_x^mh\bigr]\,dx,
\end{align*}
which for $m=1,2$ equals the nonnegative integrals
\begin{align*}
\int_{\R}\omega W_x(2\phi_x^2+h_x^2)\,dx \qquad\text{and}\qquad \int_{\R}\omega W_x\bigl[\phi_{xx}^2+(\phi_{xx}-h_{xx})^2\bigr]\,dx.
\end{align*}
Replacing the frozen coefficients by $(V+\phi,h)$ produces an error bounded by $CA\int_\R W_x(\phi_x^2+h_x^2+\phi_{xx}^2+h_{xx}^2)\,dx$. In $\kappa$ times \eqref{eqs5.8}, the term $\kappa\int_\R\omega\mu(V+\phi)h_{xxx}\phi_{xx}\,dx$ is bounded by $\frac{\kappa}{4}\int_\R\omega\frac{w(V+\phi)^2}{W}\phi_{xx}^2\,dx+C\kappa\int_\R\omega\mu(V+\phi)h_{xxx}^2\,dx$, and $\kappa\int_\R \omega W h_{xx}^2\,dx$ is absorbed by $\int_\R \omega\mu(V+\phi)h_{xx}^2\,dx$ for small $\kappa>0$.

\medskip
\noindent\textbf{Step 4.}
Expanding the derivatives of $\mu(V+\phi)$ and $w(V+\phi)^2/W$ by the chain rule, we note that in $\p_t(w(V+\phi)^2/W^2)$ the coefficient multiplying $W_t$ is $O(|\phi|)$, while the $R_1$ term satisfies
\begin{align*}
A^2\int_{\R}|R_1|(\phi_x^2+\phi_{xx}^2)\,dx\le C\|R_1\|_{H^1}A^2E_1(t)\le C\mathcal{B}(t)A^2E_1(t).
\end{align*}
Every nonlinear term in \eqref{eqs5.7} and \eqref{eqs5.8} free of background derivatives ($W_x,W_{xx},\dots$) and residuals ($R_1,R_2$) is bounded by a constant multiple of one of the nine products
\begin{align*}
&|h_x|\phi_x^2,\quad |\phi_xh_xh_{xx}|,\quad (|\phi_x|+|h_x|)(\phi_{xx}^2+h_{xx}^2+h_{xxx}^2), \quad |\phi_x|^3|h_{xx}|, \\
&\phi_x^2|h_xh_{xxx}|,\quad \phi_x^2|\phi_{xx}|, \quad \phi_x^2|h_x\phi_{xx}|, \quad |\phi_xh_{xx}\phi_{xx}|, \quad |h_x|\phi_{xx}^2.
\end{align*}

To bound $|h_x|\phi_x^2$ and $|\phi_xh_xh_{xx}|$, we relate $\int_\R (|\phi_x|+|h_x|)^3\,dx$ to $\int_\R z_\pm^2(z_\pm)_x^2\,dx$. Differentiating \eqref{eqs3.1} once and twice in $x$ gives
\begin{align}\label{eqs5.10}
\begin{cases}
(z_-)_x=Wh_x+W_xh+w(V+\phi)\phi_x+[w(V+\phi)-W]V_x,\\
(z_+)_x=Wh_x+W_xh-w(V+\phi)\phi_x-[w(V+\phi)-W]V_x,
\end{cases}
\end{align}
and
\begin{align}\label{eqs5.11}
\begin{cases}
(z_-)_{xx}=Wh_{xx}+2W_xh_x+W_{xx}h+w(V+\phi)\phi_{xx}+w'(V+\phi)\phi_x^2+2w'(V+\phi)V_x\phi_x\\
\qquad\qquad+[w'(V+\phi)-w'(V)]V_x^2+[w(V+\phi)-W]V_{xx},\\[4pt]
(z_+)_{xx}=Wh_{xx}+2W_xh_x+W_{xx}h-w(V+\phi)\phi_{xx}-w'(V+\phi)\phi_x^2-2w'(V+\phi)V_x\phi_x\\
\qquad\qquad-[w'(V+\phi)-w'(V)]V_x^2-[w(V+\phi)-W]V_{xx}.
\end{cases}
\end{align}
For $z \in \{z_-, z_+\}$, integration by parts and the Cauchy--Schwarz inequality yield
\begin{align*}
\int_{\mathbb{R}}\left|z_x\right|^3 d x=\int_{\mathbb{R}}\left|z_x\right| z_x z_x d x=-\int_{\mathbb{R}}\left(\left|z_x\right| z_x\right)_x z d x=-2 \int_{\mathbb{R}} z\left|z_x\right| z_{x x} d x \leq 2\left(\int_{\mathbb{R}} z^2 z_x^2 d x\right)^{1 / 2}\left\|z_{x x}\right\|_{L^2}.
\end{align*}
From \eqref{eqs5.10}--\eqref{eqs5.11}, $|V_x|\le CW_x$, and $|V_{xx}|\le C(|W_{xx}|+W_x^2)$, we have
\[
|\phi_x|+|h_x|\le C\bigl(|(z_-)_x|+|(z_+)_x|\bigr)+CW_x(|\phi|+|h|)
\]
and
\begin{align}\label{eqs5.12}
\|(z_-)_{xx}\|_{L^2}+\|(z_+)_{xx}\|_{L^2}&\le C\bigl[\|\phi_{xx}\|_{L^2}+\|h_{xx}\|_{L^2}+\|\phi_x^2\|_{L^2}\nonumber\\
&\quad+\|W_x\|_{L^\infty}\|(\phi_x,h_x)\|_{L^2}+A\bigl(\|W_{xx}\|_{L^2}+\|(W_x)^2\|_{L^2}\bigr)\bigr].
\end{align}
Since $Q=\mu(V+\phi)\phi_x+[\mu(V+\phi)-\mu(V)]V_x$, \eqref{eqs3.8} gives
\begin{align}\label{eqs5.13}
\int_{\R}\phi_x^4\,dx\le C\int_{\R} Q^4\,dx+C\|W_x\|_{L^\infty}^4E(t).
\end{align}
Moreover, by \eqref{eqs2.10}, \eqref{eqs2.20}, and $\delta\le A^2\le 1$,
\begin{align*}
A^2\int_{\R} W_x^3(|\phi|+|h|)^3\,dx\le CA^5\|W_x\|_{L^\infty}^2\int_{\R} W_x\,dx\le CA^5\delta(1+t)^{-2}\le CA^4(1+t)^{-3/2}.
\end{align*}
Combining these bounds and applying Young's inequality, we find that for any $\beta>0$ and sufficiently small $A$,
\begin{align}\label{eqs5.14}
A^2\int_{\R}(|\phi_x|+|h_x|)^3\,dx&\le \beta\int_{\R}\bigl[z_-^2(z_-)_x^2+z_+^2(z_+)_x^2+Q^4\bigr]\,dx+\frac{c}{16}A^2D_1(t)\nonumber\\
&\quad+C_\beta\mathcal B(t)\bigl(E(t)+A^2E_1(t)\bigr)+CA^4(1+t)^{-3/2}.
\end{align}
Using $\|\phi_x\|_{L^\infty}+\|h_x\|_{L^\infty}\le CA$ from \eqref{eqs2.20} together with \eqref{eqs5.13}--\eqref{eqs5.14}, we bound $A^2$ times the integrals of the nine products above by
\begin{align*}
&A^2\int_{\R}\bigl[|h_x|\phi_x^2+|\phi_x h_x h_{xx}|+(|\phi_x|+|h_x|)(\phi_{xx}^2+h_{xx}^2+h_{xxx}^2)\bigr]\,dx\\
&\quad\le CA^2\int_{\R}(|\phi_x|+|h_x|)^3\,dx+\beta A^2\|h_{xx}\|_{L^2}^2+CA^3\|(\phi_{xx},h_{xx},h_{xxx})\|_{L^2}^2,
\end{align*}
and
\begin{align*}
&A^2\int_{\R}\bigl(|\phi_x|^3|h_{xx}|+\phi_x^2|h_x h_{xxx}|+\phi_x^2|\phi_{xx}|+\phi_x^2|h_x\phi_{xx}|+|\phi_x h_{xx}\phi_{xx}|+|h_x|\phi_{xx}^2\bigr)\,dx\\
&\quad\le \beta A^2\|(\phi_{xx},h_{xx},h_{xxx})\|_{L^2}^2+CA^3\|(\phi_{xx},h_{xx})\|_{L^2}^2+C_\beta\bigl(A^2+A^4\bigr)\int_{\R}\phi_x^4\,dx,
\end{align*}
all of which are controlled by the right-hand side of \eqref{eqs5.6} for small $\beta$ and $A$.

\medskip
\noindent\textbf{Step 5.}
Using the relations
\begin{align}\label{eqs5.15}
V_x=\frac{W_x}{w'(V)}, \qquad V_{xx}=\frac{W_{xx}}{w'(V)}-\frac{w''(V)}{w'(V)^3}W_x^2,\qquad \p_t\omega=4\omega\Bigl(\frac{W_{xx}}{W}-W_x\Bigr),\qquad \p_x\omega=\frac{4\omega W_x}{W},
\end{align}
we group the remaining profile and residual terms in the evolution of $A^2E_1(t)$ into three sums $\mathcal{J}_1, \mathcal{J}_2, \mathcal{J}_3$:

\smallskip
\noindent(1) Let $\mathcal{J}_1$ be the sum of terms where profile coefficients multiply differentiated perturbations. After moving the principal $W_x$-quadratic forms to the left-hand side as in Step~3, the integrand of $\mathcal{J}_1$ is bounded by
\[
\begin{aligned}
&CA^2\bigl(W_x+|W_{xx}|+|W_{xxx}|+W_x^2+|R_1|\bigr)\bigl(|\phi_x|+|h_x|+|\phi_{xx}|+|h_{xx}|\bigr)\\
&\quad\times\bigl(|\phi_x|+|h_x|+|\phi_{xx}|+|h_{xx}|+|h_{xxx}|\bigr),
\end{aligned}
\]
so Cauchy--Schwarz, Young's inequality, and \eqref{eqs5.4} give
\begin{equation}\label{eqs5.16}
|\mathcal{J}_1|\le \frac{\beta}{4}A^2D_1(t)+C_\beta\mathcal B(t)A^2E_1(t).
\end{equation}

The above bound of $\mathcal{J}_1$ is written in a unified form for convenience.
In the actual expansion, after extracting the principal quadratic
forms in Step~3, the single-$W_x$ terms involving only $\phi_x$
and $h_x$ carry an additional factor $O(A)$; in particular,
no term of the form $CA^2W_x(\phi_x^2+h_x^2)$  remains.
The remaining single-\(W_x\) terms contain at least one of \(\phi_{xx}\), \(h_{xx}\), or \(h_{xxx}\), and can be controlled as in \eqref{eqs5.16}.

\smallskip
\noindent(2) Let $\mathcal{J}_2$ be the sum of terms where profile derivatives multiply $(\phi,h)$, which are bounded by integrals of
\begin{align*}
A^2\bigl(W_x^2+|W_{xx}|+|W_{xxx}|+|W_{xxxx}|\bigr)(|\phi|+|h|)\bigl(|\phi_x|+|h_x|+|\phi_{xx}|+|h_{xx}|+|h_{xxx}|\bigr).
\end{align*}
When $W_{xx}$ multiplies $(\phi_{xx},h_{xx})$ and $(\phi,h)$, we integrate by parts once:
\begin{align*}
\int_{\R}\omega W_{xx}\phi\phi_{xx}\,dx=-\int_{\R}\omega W_{xxx}\phi\phi_x\,dx-\int_{\R}\omega W_{xx}\phi_x^2\,dx-4\int_{\R}\omega\frac{W_x}{W}W_{xx}\phi\phi_x\,dx,
\end{align*}
and similarly for $h$. Using $\|(\phi,h)\|_{L^\infty}\le CA$, the interpolation $\|(\phi_x,h_x)\|_{L^2}\le CA^{1/2}\|(\phi_{xx},h_{xx})\|_{L^2}^{1/2}$, and Young's inequality, we obtain
\begin{align}\label{eqs5.17}
|\mathcal{J}_2|&\le \frac{\beta}{4}A^2D_1(t)+C_\beta\mathcal B(t)A^2E_1(t)\nonumber\\
&\quad+C_\beta A^4\Bigl(\|W_{xx}\|_{H^2}^2+\|(W_x)^2\|_{H^2}^2+\|(W_x)^2\|_{L^2}^{4/3}+\|W_{xx}\|_{L^2}^{4/3}+\|W_{xxx}\|_{L^2}^{4/3}\Bigr).
\end{align}

\smallskip
\noindent(3) Let $\mathcal{J}_3$ denote the residual terms in \eqref{eqs5.7} and \eqref{eqs5.8} involving $R_1$ and $R_2$. Integrating the first-order pairings $\int_\R \omega\frac{w^2}{W^2}\phi_x(R_1)_x\,dx$ and $\int_\R \omega h_x(R_2/W)_x\,dx$ by parts once and applying Cauchy--Schwarz yields
\begin{equation}\label{eqs5.18}
|\mathcal{J}_3|\le \frac{\beta}{4}A^2D_1(t)+C_\beta\mathcal B(t)A^2E_1(t)+C_\beta A^2\|(R_1,R_2)\|_{H^2}^2.
\end{equation}

By Lemma~\ref{lem2.1} and $\delta\le A^2$, the profile and residual norms in \eqref{eqs5.17}--\eqref{eqs5.18} are bounded by $CA^4(1+t)^{-3/2}$. Combining \eqref{eqs5.14} and \eqref{eqs5.16}--\eqref{eqs5.18} gives \eqref{eqs5.6}.
\end{proof}
At the end of this section, we introduce the one-dimensional Nash inequality.
\begin{Lem}[One-dimensional Nash inequality \cite{Nash1958}]\label{lem5.1add}
Let $u\in H^1(\R)\cap L^p(\R)$ with $p\ge4$. Then
\begin{align}\label{add5.1}
\Bigl(\int_\R|u|^p\,dx\Bigr)^{\frac{p+2}{p-2}}\le C_p\Bigl(\int_\R|u|^2\,dx\Bigr)^{\frac{2p}{p-2}}\int_\R|u|^{p-2}u_x^2\,dx.
\end{align}
\end{Lem}
\noindent Moreover, the derivative energy estimate is compatible with the Nash inequality as follows:
\begin{Lem}\label{lem5.2}
By Theorem~\ref{thm2.1}, we have
\begin{align}\label{eqs5.19}
\|(\phi,h)(t)\|_{H^2}+\|Q(t)\|_{H^1}+\|z_-(t)\|_{L^2}+\|z_+(t)\|_{L^2}\le CA,\qquad 0<A\le 1,
\end{align}
\begin{align}\label{eqs5.20}
E(t)+A^2E_1(t)\simeq\|z_-(t)\|_{L^4}^4+\|z_+(t)\|_{L^4}^4+\|Q(t)\|_{L^4}^4+A^2\|(\phi_x,h_x)(t)\|_{H^1}^2,
\end{align}
and
\begin{equation}\label{eqs5.21}
\bigl(E(t)+A^2E_1(t)\bigr)^3\le CA^8\Bigl(\int_{\R}\bigl[z_-^2(z_-)_x^2+z_+^2(z_+)_x^2+Q^4\bigr]\,dx+A^2D_1(t)\Bigr).
\end{equation}
\end{Lem}

\begin{proof}
Differentiating $Q=\mu(V+\phi)\phi_x+[\mu(V+\phi)-\mu(V)]V_x$ in $x$ and using \eqref{eqs2.9}--\eqref{eqs2.10} and Theorem~\ref{thm2.1} gives $\|Q\|_{H^1}\le C\|(\phi,\psi)\|_{H^2}+C\delta\le CA$, proving \eqref{eqs5.19}. The equivalence \eqref{eqs5.20} follows from \eqref{eqs5.5} and $|z-Q|^4+|z|^4\simeq |z|^4+|Q|^4$.

Applying Lemma~\ref{lem5.1add} with $p=4$ to $u=z_\pm$ and using $\|z_\pm\|_{L^2}^2\le CA^2$ gives
\begin{align}\label{eqs5.22}
\Bigl(\int_{\R} z_\pm^4\,dx\Bigr)^3\le C\Bigl(\int_{\R} z_\pm^2\,dx\Bigr)^4\int_{\R} z_\pm^2(z_\pm)_x^2\,dx\le CA^8\int_{\R}\bigl[z_-^2(z_-)_x^2+z_+^2(z_+)_x^2+Q^4\bigr]\,dx.
\end{align}
By Sobolev embedding, $\int_\R Q^4\,dx\le \|Q\|_{L^\infty}^2\|Q\|_{L^2}^2\le CA^4$, so
\begin{align}\label{eqs5.23}
\Bigl(\int_{\R} Q^4\,dx\Bigr)^3\le CA^8\int_{\R} Q^4\,dx\le CA^8\int_{\R}\bigl[z_-^2(z_-)_x^2+z_+^2(z_+)_x^2+Q^4\bigr]\,dx.
\end{align}
Finally, integration by parts and $\|(\phi,h)\|_{L^2}\le CA$ give
\begin{equation}\label{eqs5.24}
\|\phi_x\|_{L^2}^4\le CA^2\|\phi_{xx}\|_{L^2}^2,
\qquad \|h_x\|_{L^2}^4\le CA^2\|h_{xx}\|_{L^2}^2.
\end{equation}
Combining \eqref{eqs5.24} with $\|(\phi_x,h_x)\|_{H^1}^2\le CA^2$ and \eqref{eqs5.5} yields
\begin{equation}\label{eqs5.25}
E_1(t)^3\le CA^4D_1(t).
\end{equation}
Multiplying \eqref{eqs5.25} by $A^6$ and combining with \eqref{eqs5.22}--\eqref{eqs5.23} yields \eqref{eqs5.21}.
\end{proof}

\section{Corrected characteristic modes and half-space identity}\label{sec6}
\subsection{Corrected characteristic modes}\label{sec6.1}
Recall the state-space parameterization $2\int_V^{v(\xi;V)}w(y)\,dy=\xi$ from Step~4 of Proposition~\ref{prop4.1}, and set $\widehat w(\xi,V):=w(v(\xi;V))$, so that $\widehat w(0,V)=W$ and $\xi=z_--z_+$ along the solution. With
\begin{align}\label{eqs6.1}
\alpha(V):=\frac{\p_\xi\widehat w(0,V)}{2W}=\frac{w'(V)}{4W^2},
\end{align}
we define for $m\in\{2,4\}$ the corrected characteristic densities $\eta_\pm^{(m)}$ and fluxes $f_\pm^{(m)}$ by
\begin{align}\label{eqs6.2}
\eta_-^{(m)}&:=z_-^m(1+\alpha(V) z_+),\qquad \eta_+^{(m)}:=z_+^m(1-\alpha(V) z_-),\nonumber\\
f_-^{(m)}&:=-m\int_0^{z_-}\widehat w(r,V)r^{m-1}\,dr+\alpha(V) z_-^m\int_0^{z_+}\widehat w(z_--r,V)\,dr,\nonumber\\
f_+^{(m)}&:=m\int_0^{z_+}\widehat w(-r,V)r^{m-1}\,dr+\alpha(V) z_+^m\int_0^{z_-}\widehat w(r-z_+,V)\,dr.
\end{align}
For $A_0$ small enough that $|\alpha(V)|(|z_-|+|z_+|)\le 1/2$, we have $\frac12z_\pm^m\le \eta_\pm^{(m)}\le \frac32z_\pm^m$ for $m\in\{2,4\}$.

By \eqref{eqs6.2}, the transverse derivatives satisfy $(f_\mp^{(m)})_{z_\pm}=\pm\widehat w(z_--z_+,V)(\eta_\mp^{(m)})_{z_\pm}$. The compatibility defects in the primary directions are
\begin{align*}
d_-^{(m)}:=(f_-^{(m)})_{z_-}+\widehat w(z_--z_+,V)(\eta_-^{(m)})_{z_-},\qquad d_+^{(m)}:=(f_+^{(m)})_{z_+}-\widehat w(z_--z_+,V)(\eta_+^{(m)})_{z_+}.
\end{align*}
Note that $d_-^{(m)}=0$ on $z_+=0$ and $d_+^{(m)}=0$ on $z_-=0$, and
\begin{align*}
\p_{z_+}d_-^{(m)}&=z_-^{m-1}\bigl[m(2\alpha\widehat w-\widehat w_{\xi})+\alpha\widehat w_{\xi}(z_--mz_+)\bigr],\\
\p_{z_-}d_+^{(m)}&=-z_+^{m-1}\bigl[m(\widehat w_{\xi}-2\alpha\widehat w)+\alpha\widehat w_{\xi}(z_+-mz_-)\bigr].
\end{align*}
By \eqref{eqs6.1}, $2\alpha(V)\widehat w(0,V)-\widehat w_\xi(0,V)=0$, so $|2\alpha\widehat w-\widehat w_\xi|\le C(|z_-|+|z_+|)$. Integrating from $0$ to $z_\pm$ gives
\begin{align}\label{eqs6.3}
|d_-^{(m)}|\le C|z_-|^{m-1}|z_+|(|z_-|+|z_+|),\qquad |d_+^{(m)}|\le C|z_+|^{m-1}|z_-|(|z_-|+|z_+|).
\end{align}

\subsection{Half-space interaction identity}\label{sec6.2}
Evaluating $\eta_\pm^{(m)}$ and $f_\pm^{(m)}$ along $(z_-(t,x),z_+(t,x);V(t,x))$, we define
\begin{equation}\label{eqs6.4}
\mathcal I(t):=\int_{\R}\eta_-^{(4)}(t,x)\Bigl(\int_{-\infty}^x\eta_+^{(2)}(t,y)\,dy\Bigr)\,dx+\int_{\R}\eta_-^{(2)}(t,x)\Bigl(\int_{-\infty}^x\eta_+^{(4)}(t,y)\,dy\Bigr)\,dx.
\end{equation}
By the chain rule,
\begin{align*}
(\eta_\pm^{(m)})_t+(f_\pm^{(m)})_x&=(\eta_\pm^{(m)})_{z_-}\bigl[(z_-)_t-w(v)(z_-)_x\bigr]+(\eta_\pm^{(m)})_{z_+}\bigl[(z_+)_t+w(v)(z_+)_x\bigr]\\
&\quad+(\eta_\pm^{(m)})_VV_t+(f_\pm^{(m)})_VV_x+d_\pm^{(m)}(z_\pm)_x.
\end{align*}
Setting $S_\pm^{(m)}:=(\eta_\pm^{(m)})_{z_-}+(\eta_\pm^{(m)})_{z_+}$, substituting \eqref{eqs3.3}--\eqref{eqs3.4}, and using $V_t = R_1-WV_x$, we obtain
\begin{align}\label{eqs6.5}
(\eta_\pm^{(m)})_t+(f_\pm^{(m)})_x&=\p_x\Bigl(S_\pm^{(m)}[\mu(v)u_x-\mu(V)U_x]\Bigr)-[\mu(v)u_x-\mu(V)U_x]\p_x S_\pm^{(m)}-R_2 S_\pm^{(m)}\nonumber\\
&\quad+R_1\bigl[-W(\eta_\pm^{(m)})_{z_-}+W(\eta_\pm^{(m)})_{z_+}+(\eta_\pm^{(m)})_V\bigr]\nonumber\\
&\quad+V_x\bigl[(f_\pm^{(m)})_{V}-W(\eta_\pm^{(m)})_{V}+2W(w(v)-W)(\eta_\pm^{(m)})_{z_+}\bigr]+d_\pm^{(m)}(z_\pm)_x\nonumber\\
&=:\mathcal{F}_\pm^{(m)}.
\end{align}
Applying the half-space identity for one-dimensional balance laws $a_t+p_x=F$, $b_t+q_x=G$,
\begin{align*}
&\frac{d}{dt}\int_{\R} a(x)\Bigl(\int_{-\infty}^xb(y)\,dy\Bigr)\,dx+\int_{\R}\bigl[a(x)q(x)-p(x)b(x)\bigr]\,dx\\
&=\int_{\R} F(x)\Bigl(\int_{-\infty}^xb(y)\,dy\Bigr)\,dx+\int_{\R} G(x)\Bigl(\int_x^\infty a(y)\,dy\Bigr)\,dx,
\end{align*}
to $(a,b)=(\eta_-^{(4)},\eta_+^{(2)})$ and $(a,b)=(\eta_-^{(2)},\eta_+^{(4)})$ and summing the two equations yields
\begin{align}\label{eqs6.6}
&\frac{d}{dt}\mathcal I(t)+\int_{\R}\bigl(\eta_-^{(4)}f_+^{(2)}-f_-^{(4)}\eta_+^{(2)}+\eta_-^{(2)}f_+^{(4)}-f_-^{(2)}\eta_+^{(4)}\bigr)\,dx\nonumber\\
&=\sum_{(m,n)\in\{(4,2),(2,4)\}}\Biggl[\int_{\R} \mathcal{F}_-^{(m)}(x)\Bigl(\int_{-\infty}^x\eta_+^{(n)}(y)\,dy\Bigr)\,dx+\int_{\R} \mathcal{F}_+^{(n)}(x)\Bigl(\int_x^{+\infty}\eta_-^{(m)}(y)\,dy\Bigr)\,dx\Biggr].
\end{align}

\begin{Lem}\label{lem6.1}
For $m\in\{2,4\}$, the quotients
\begin{align*}
\frac{\eta_-^{(m)}}{z_-^m}=1+\alpha(V) z_+,\qquad \frac{\eta_+^{(m)}}{z_+^m}=1-\alpha(V) z_-,
\end{align*}
and
\begin{align*}
\frac{f_-^{(m)}}{z_-^m}&=-m\int_0^{1}\widehat w(rz_-,V)r^{m-1}\,dr+\alpha(V) \int_0^{z_+}\widehat w(z_--r,V)\,dr,\\
\frac{f_+^{(m)}}{z_+^m}&=m\int_0^{1}\widehat w(-rz_+,V)r^{m-1}\,dr+\alpha(V)\int_0^{z_-}\widehat w(r-z_+,V)\,dr
\end{align*}
have uniformly bounded state derivatives up to second order on the small state ball, and satisfy
\begin{equation}\label{eqs6.7}
\Bigl|\frac{\eta_-^{(m)}}{z_-^m}-1\Bigr|+\Bigl|\frac{f_-^{(m)}}{z_-^m}+W\Bigr|\le C(|z_-|+|z_+|),\qquad \Bigl|\frac{\eta_+^{(m)}}{z_+^m}-1\Bigr|+\Bigl|\frac{f_+^{(m)}}{z_+^m}-W\Bigr|\le C(|z_-|+|z_+|).
\end{equation}
\end{Lem}

\begin{proof}
Changing variables $r\mapsto rz_\pm$ in \eqref{eqs6.2} gives the stated formulas. Since $\widehat w(0,V)=W$ and $m\int_0^1 r^{m-1}\,dr=1$, the mean-value theorem yields \eqref{eqs6.7}.
\end{proof}

From $\frac12 z_\pm^m\le \eta_\pm^{(m)}\le \frac32 z_\pm^m$ and $\|z_\pm(t)\|_{L^2}^2\le CA^2$, we have
\begin{align}\label{eqs6.8}
0\le\mathcal I(t)\le CA^2E(t),
\end{align}
and by \eqref{eqs6.7}, for sufficiently small $A$,
\begin{align}\label{eqs6.9}
\eta_-^{(4)}f_+^{(2)}-f_-^{(4)}\eta_+^{(2)}+\eta_-^{(2)}f_+^{(4)}-f_-^{(2)}\eta_+^{(4)}&=2W(z_-^4z_+^2+z_-^2z_+^4)+O\bigl((|z_-|+|z_+|)(z_-^4z_+^2+z_-^2z_+^4)\bigr)\nonumber\\
&\ge c(z_-^4z_+^2+z_-^2z_+^4).
\end{align}

\section{Nonlocal wave-interaction estimates}\label{sec7}

\begin{Prop}\label{prop7.1}
For every $\beta>0$ and sufficiently small $A,\delta$, the functional $\mathcal{I}(t)$ satisfies
\begin{align}\label{eqs7.1}
&\frac{d}{dt}\mathcal I(t)+c\int_{\R} (z_-^4z_+^2+z_-^2z_+^4)\,dx \notag \\
&\le \beta\int_{\R}\bigl[z_-^2(z_-)_x^2+z_+^2(z_+)_x^2+Q^4\bigr]\,dx+\beta\int_{\R} W_x(z_--z_+)^2\bigl[(z_-+z_+)^2+2z_+^2\bigr]\,dx \notag \\
&\quad+(\beta+C_\beta A^2)\int_{\R}(z_-^4z_+^2+z_-^2z_+^4)\,dx+C_\beta A^2D_1(t)+C_\beta\mathcal B(t)\bigl(E(t)+A^2E_1(t)\bigr) \notag \\
&\quad+CA\bigl(\|R_1\|_{L^2}+\|R_2\|_{L^2}\bigr)E(t)+CA^4(1+t)^{-3/2}+CA^2\bigl(\|R_1\|_{L^4}+\|R_2\|_{L^4}\bigr)E(t)^{3/4}.
\end{align}
\end{Prop}

\begin{proof}
We substitute \eqref{eqs6.5} into \eqref{eqs6.6} and estimate each contribution in turn.

\medskip
\noindent\textbf{Step 1.}
For the quartic defects paired against the quadratic tails $\int_{-\infty}^x \eta_+^{(2)}\,dy\le CA^2$ and $\int_x^\infty \eta_-^{(2)}\,dy\le CA^2$, \eqref{eqs6.3}, $\|z_\pm\|_{L^\infty}\le CA$, and Young's inequality give
\begin{align*}
&\Bigl|\int_{\R} d_-^{(4)}(z_-)_x\Bigl(\int_{-\infty}^x\eta_+^{(2)}(y)\,dy\Bigr)\,dx\Bigr|+\Bigl|\int_{\R} d_+^{(4)}(z_+)_x\Bigl(\int_x^\infty\eta_-^{(2)}(y)\,dy\Bigr)\,dx\Bigr|\\
&\le CA^2\int_{\R} |z_-|^3|z_+|(|z_-|+|z_+|)|(z_-)_x|\,dx+CA^2\int_{\R} |z_+|^3|z_-|(|z_-|+|z_+|)|(z_+)_x|\,dx\\
&\le \beta\int_{\R}\bigl[z_-^2(z_-)_x^2+z_+^2(z_+)_x^2+Q^4\bigr]\,dx+C_\beta A^6\int_{\R}(z_-^4z_+^2+z_-^2z_+^4)\,dx.
\end{align*}
For the quadratic defects paired against the quartic tails $\int_{\R} \eta_\pm^{(4)}\,dy\le C\|z_\pm\|_{L^4}^4$, \eqref{eqs6.3}, Young's inequality, and the Nash inequality \eqref{eqs5.22} give
\begin{align*}
&\Bigl|\int_{\R} d_-^{(2)}(z_-)_x\Bigl(\int_{-\infty}^x\eta_+^{(4)}(y)\,dy\Bigr)\,dx\Bigr|+\Bigl|\int_{\R} d_+^{(2)}(z_+)_x\Bigl(\int_x^\infty\eta_-^{(4)}(y)\,dy\Bigr)\,dx\Bigr|\\
&\le C\|z_+\|_{L^4}^4\int_{\R} |z_-||z_+|(|z_-|+|z_+|)|(z_-)_x|\,dx+C\|z_-\|_{L^4}^4\int_{\R} |z_+||z_-|(|z_-|+|z_+|)|(z_+)_x|\,dx\\
&\le \beta\int_{\R}\bigl[z_-^2(z_-)_x^2+z_+^2(z_+)_x^2\bigr]\,dx+C_\beta\bigl(\|z_-\|_{L^4}^4+\|z_+\|_{L^4}^4\bigr)^3\\
&\le (\beta+C_\beta A^8)\int_{\R}\bigl[z_-^2(z_-)_x^2+z_+^2(z_+)_x^2+Q^4\bigr]\,dx.
\end{align*}

\medskip
\noindent\textbf{Step 2.}
Integrating the divergence term in \eqref{eqs6.5} by parts against the cumulative tails in \eqref{eqs6.6} produces $[\mu(v)u_x-\mu(V)U_x]$ multiplied by
\[
\sum_{(m,n)\in\{(4,2),(2,4)\}}
 \bigl(\eta_-^{(m)}S_+^{(n)}-\eta_+^{(n)}S_-^{(m)}\bigr).
\]
By \eqref{eqs4.2} and $V_x=W_x/w'(V)$,
\begin{align}\label{eqs7.2}
2[\mu(v)u_x-\mu(V)U_x]=\mu(v)\bigl[(z_-)_x+(z_+)_x\bigr]+\frac{2W[\mu(V)-\mu(v)]}{w'(V)}W_x.
\end{align}
By Lemma~\ref{lem6.1}, the four products $\eta_-^{(4)}S_+^{(2)}$, $\eta_+^{(2)}S_-^{(4)}$, $\eta_-^{(2)}S_+^{(4)}$, and $\eta_+^{(4)}S_-^{(2)}$ are bounded by
\begin{align*}
C|z_-|^4|z_+|,\qquad C|z_-|^3z_+^2,\qquad Cz_-^2|z_+|^3,\qquad C|z_-|z_+^4,
\end{align*}
respectively. When multiplied by $\frac12\mu(v)[(z_-)_x+(z_+)_x]$, six of the resulting eight products (whose integral we denote by $I_{\mathrm{loc},1}$) carry at least two powers of the differentiated variable, namely $|z_-|^4|z_+||(z_-)_x|$, $|z_-|^3z_+^2|(z_\pm)_x|$, $z_-^2|z_+|^3|(z_\pm)_x|$, and $|z_-|z_+^4|(z_+)_x$. By Young's inequality and $\|z_\pm\|_{L^\infty}\le CA$,
\begin{align}\label{eqs7.3}
|I_{\mathrm{loc},1}|
\le \beta\int_{\R}\bigl[z_-^2(z_-)_x^2+z_+^2(z_+)_x^2\bigr]\,dx+C_\beta A^2\int_{\R}(z_-^4z_+^2+z_-^2z_+^4)\,dx.
\end{align}
For the remaining two endpoint terms $I_{\mathrm{loc},2}$, where $\eta_-^{(4)}S_+^{(2)}$ multiplies $(z_+)_x$ and $-\eta_+^{(4)}S_-^{(2)}$ multiplies $(z_-)_x$, we integrate in state space to transfer the spatial derivative onto the fourth-order factor:
\begin{align*}
-\frac12\int_{\R}\mu(v)\eta_+^{(4)}S_-^{(2)}(z_-)_x\,dx
&=\frac12\int_{\R}\p_{z_+}\Bigl[\int_0^{z_-}\mu(v(r-z_+;V))\eta_+^{(4)}(r,z_+;V)S_-^{(2)}(r,z_+;V)\,dr\Bigr](z_+)_x\,dx\\
&\quad+\frac12\int_{\R}\p_{V}\Bigl[\int_0^{z_-}\mu(v(r-z_+;V))\eta_+^{(4)}(r,z_+;V)S_-^{(2)}(r,z_+;V)\,dr\Bigr]V_x\,dx.
\end{align*}
Since $\eta_+^{(4)}(r,z_+;V)=O(z_+^4)$ and $S_-^{(2)}(r,z_+;V)=O(|r|)$, the two brackets above are bounded by $C|z_-|^2|z_+|^3$ and $C|z_-|^2|z_+|^4$. Treating $\frac12\int_\R \mu(v)\eta_-^{(4)}S_+^{(2)}(z_+)_x\,dx$ symmetrically by integrating in $z_+$, we obtain
\begin{align}\label{eqs7.4}
|I_{\mathrm{loc},2}|
\le \beta\int_{\R}\bigl[z_-^2(z_-)_x^2+z_+^2(z_+)_x^2\bigr]\,dx+(C_\beta A^2+C\|W_x\|_{L^\infty})\int_{\R}(z_-^4z_+^2+z_-^2z_+^4)\,dx.
\end{align}
For the $W_x$-term in \eqref{eqs7.2}, since $|\mu(V)-\mu(v)|\le C|z_--z_+|$, its product with the four mode-derivative terms is bounded pointwise by
\begin{align*}
CW_x|z_--z_+||z_-z_+|(|z_-|+|z_+|)^3\le\beta W_x(z_--z_+)^2\bigl[(z_-+z_+)^2+2z_+^2\bigr]+C_\beta A^2(z_-^4z_+^2+z_-^2z_+^4).
\end{align*}

\medskip
\noindent\textbf{Step 3.}
Next consider the terms in \eqref{eqs6.6} coming from $-[\mu(v)u_x-\mu(V)U_x]\p_x S_\pm^{(m)}$ multiplied by the cumulative tails. For $m\in\{2,4\}$, differentiating \eqref{eqs6.2} gives
\begin{align}\label{eqs7.5}
|\p_x S_-^{(m)}|&\le C\bigl[|z_-|^{m-2}|(z_-)_x|+|z_-|^{m-1}|(z_+)_x|+W_x|z_-|^{m-1}(|z_-|+|z_+|)\bigr],\nonumber\\
|\p_x S_+^{(m)}|&\le C\bigl[|z_+|^{m-2}|(z_+)_x|+|z_+|^{m-1}|(z_-)_x|+W_x|z_+|^{m-1}(|z_-|+|z_+|)\bigr],
\end{align}
while the cumulative tails satisfy
\begin{align*}
\Bigl\|\int_{-\infty}^x\eta_+^{(2)}(y)\,dy\Bigr\|_{L^\infty}+\Bigl\|\int_x^\infty\eta_-^{(2)}(y)\,dy\Bigr\|_{L^\infty}&\le CA^2,\\
\Bigl\|\int_{-\infty}^x\eta_+^{(4)}(y)\,dy\Bigr\|_{L^\infty}+\Bigl\|\int_x^{\infty}\eta_-^{(4)}(y)\,dy\Bigr\|_{L^\infty}&\le CE(t).
\end{align*}
When $m=4$, multiplying \eqref{eqs7.2} by \eqref{eqs7.5} against the quadratic tails of size $CA^2$ produces $A^2\int_\R [z_-^2(z_-)_x^2+z_+^2(z_+)_x^2]\,dx$ together with the cross-family gradient terms
\[
A^2\int_\R \bigl[z_-^2(z_+)_x^2+z_+^2(z_-)_x^2+|z_-|^3(z_+)_x^2+|z_+|^3(z_-)_x^2\bigr]\,dx.
\]
Integrating by parts in $x$ (e.g., $\int_\R z_-^2(z_+)_x^2\,dx=-2\int_\R z_-(z_-)_x z_+(z_+)_x\,dx-\int_\R z_-^2 z_+(z_+)_{xx}\,dx$) and applying Young's inequality with $\|z_\pm\|_{L^\infty}\le CA$, we obtain
\begin{align}\label{eqs7.7}
&A^2\int_{\R}\bigl[z_-^2(z_+)_x^2+z_+^2(z_-)_x^2+|z_-|^3(z_+)_x^2+|z_+|^3(z_-)_x^2\bigr]\,dx\nonumber\\
&\le \beta\int_{\R}\bigl[z_-^2(z_-)_x^2+z_+^2(z_+)_x^2\bigr]\,dx+\beta\int_{\R}(z_-^4z_+^2+z_-^2z_+^4)\,dx+C_\beta A^4\bigl(\|(z_-)_{xx}\|_{L^2}^2+\|(z_+)_{xx}\|_{L^2}^2\bigr).
\end{align}
By \eqref{eqs5.11} and \eqref{eqs5.12},
\begin{align}\label{eqs7.8}
\|(z_-)_{xx}\|_{L^2}^2+\|(z_+)_{xx}\|_{L^2}^2&\le CD_1(t)+C\bigl(\|(\phi_x,\psi_x)\|_{L^\infty}^2+\|W_x\|_{L^\infty}^2\bigr)E_1(t)\nonumber\\
&\quad+CA^2\bigl(\|W_{xx}\|_{L^2}^2+\|W_{x}^2\|_{L^2}^2+\|R_1\|_{H^1}^2\bigr),
\end{align}

\noindent which leads to 
\begin{equation}
\label{eqs_step3_m=4}
   \begin{aligned}
& A^2 \int_{\mathbb{R}}\left[z_{-}^2\left(z_{+}\right)_x^2+z_{+}^2\left(z_{-}\right)_x^2+\left|z_{-}\right|^3\left(z_{+}\right)_x^2+\left|z_{+}\right|^3\left(z_{-}\right)_x^2\right] d x \\
\leq & \beta \int_{\mathbb{R}}\left[z_{-}^2\left(z_{-}\right)_x^2+z_{+}^2\left(z_{+}\right)_x^2\right] d x+\beta \int_{\mathbb{R}}\left(z_{-}^4 z_{+}^2+z_{-}^2 z_{+}^4\right) d x+C_\beta A^4 D_1(t) \\
& +C_\beta A^4 \mathcal{B}(t) E_1(t)+C_\beta A^6\left(\left\|W_{x x}\right\|_{L^2}^2+\left\|W_x^2\right\|_{L^2}^2+\left\|R_1\right\|_{H^1}^2\right).
\end{aligned}
\end{equation}
When $m=2$, the cumulative viscous integral $I_{\mathrm{cum},2}$ is paired against the quartic tails of size $CE(t)$. Since $|\mu(v)u_x-\mu(V)U_x|\le C(|\psi_x|+W_x|\phi|)$ and $|\p_x S_\pm^{(2)}|\le C(|(z_-)_x|+|(z_+)_x|+W_x(|\phi|+|\psi|))$, the Cauchy--Schwarz inequality gives
\begin{equation}\label{eqs7.9}
|I_{\mathrm{cum},2}|
\le C\Bigl[\|\phi_x\|_{L^2}^2+\|\psi_x\|_{L^2}^2+\int_{\R} W_x(\phi^2+\psi^2)\,dx\Bigr] E(t)\le C\mathcal{B}(t)E(t).
\end{equation}
Similarly, the $W_x$-terms in \eqref{eqs7.5} for $m=4$ and the $W_x$-part of \eqref{eqs7.2} in the $m=4$ cumulative integral are bounded by
\begin{align}\label{eqs7.10}
&CA^2\bigl(\|\psi_x\|_{L^\infty}^2+\|\phi_x\|_{L^\infty}^2+\|W_x\|_{L^\infty}^2\bigr)E(t)+\beta\int_{\R}\bigl[z_-^2(z_-)_x^2+z_+^2(z_+)_x^2\bigr]\,dx\nonumber\\
&\quad+\beta\int_{\R} W_x(z_--z_+)^2\bigl[(z_-+z_+)^2+2z_+^2\bigr]\,dx.
\end{align}

\medskip
\noindent\textbf{Step 4.}
Let $I_{\mathrm{res},4}$ and $I_{\mathrm{res},2}$ denote the residual integrals of the quartic and quadratic modes, respectively. H\"older's inequality and Cauchy--Schwarz yield
\begin{align}
|I_{\mathrm{res},4}|&\le CA^2\bigl(\|R_1\|_{L^4}+\|R_2\|_{L^4}\bigr)E(t)^{3/4},\label{eqs7.11}\\
|I_{\mathrm{res},2}|&\le CA\bigl(\|R_1\|_{L^2}+\|R_2\|_{L^2}\bigr)E(t).\label{eqs7.12}
\end{align}
For the sum $I_{\mathrm{drift}}$ of the four $V_x$-drift integrals in \eqref{eqs6.6}, Lemma~\ref{lem6.1} implies
\begin{align*}
|V_x|\bigl|(f_\pm^{(2)})_V-W(\eta_\pm^{(2)})_V+2W(w(v)-W)(\eta_\pm^{(2)})_{z_+}\bigr|&\le CW_x(\phi^2+\psi^2),\\
|V_x|\bigl|(f_\pm^{(4)})_V-W(\eta_\pm^{(4)})_V+2W(w(v)-W)(\eta_\pm^{(4)})_{z_+}\bigr|&\le CW_x(z_-^4+z_+^4),
\end{align*}
so pairing with the respective tails gives
\begin{equation}\label{eqs7.13}
|I_{\mathrm{drift}}|\le C\Bigl(\int_{\R} W_x(\phi^2+\psi^2)\,dx\Bigr)E(t)+CA^2\int_{\R} W_x(z_-^4+z_+^4)\,dx.
\end{equation}
Using the pointwise bounds
\begin{align*}
z_-^4+z_+^4\le C\Bigl((z_--z_+)^2\bigl[(z_-+z_+)^2+2z_+^2\bigr]+z_-^2z_+^2\Bigr),\qquad z_-^2z_+^2\le(z_-^4z_+^2+z_-^2z_+^4)^{2/3},
\end{align*}
and applying Young's inequality with exponents $(3,3/2)$ to $A^2 W_x (z_-^4z_+^2+z_-^2z_+^4)^{2/3}$, we obtain
\begin{align*}
CA^2\int_{\R} W_x(z_-^4+z_+^4)\,dx&\le CA^2\int_{\R} W_x(z_--z_+)^2\bigl[(z_-+z_+)^2+2z_+^2\bigr]\,dx\\
&\quad+\beta\int_{\R}(z_-^4z_+^2+z_-^2z_+^4)\,dx+C_\beta A^6\int_{\R} W_x^3\,dx.
\end{align*}
Since $\int_\R W_x\,dx=\delta\le A^2$ and $\|W_x(t)\|_{L^\infty}\le C(1+t)^{-1}$, we have $A^6\int_\R W_x^3\,dx\le CA^4(1+t)^{-3/2}$, while the profile error from \eqref{eqs_step3_m=4} satisfies
\begin{align*}
C_\beta A^6\bigl(\|W_{xx}(t)\|_{H^1}^2+\|W_{x}^2(t)\|_{L^2}^2+\|R_1(t)\|_{H^1}^2\bigr)\le C_\beta A^6(1+t)^{-5/2}\min\{1,\delta\sqrt{1+t}\}^2\le C_\beta A^4(1+t)^{-3/2}.
\end{align*}
Combining \eqref{eqs6.6}, \eqref{eqs6.9}, and Steps~1--4 yields \eqref{eqs7.1}.
\end{proof}

\section{Proof of Theorem \ref{thm2.2}}\label{sec8}

\begin{proof}[Proof of Theorem~\ref{thm2.2}]
Assume without loss of generality that $0<A\le A_0$. For constants $\Lambda, M>0$, define
\begin{align}\label{eqs8.1}
E_{\mathrm{tot}}(t):=E(t)+MA^2E_1(t)+\Lambda\mathcal I(t).
\end{align}
By \eqref{eqs6.8}, for any fixed $\Lambda,M>0$ and sufficiently small $A$,
\begin{align}\label{eqs8.2}
E_{\mathrm{tot}}(t)\simeq E(t)+A^2E_1(t).
\end{align}
Set
\begin{align}\label{eqs8.3}
k(t)&:=C\mathcal B(t)+C\|W_x(t)\|_{L^\infty}\bigl(\|z_-(t)\|_{L^\infty}+\|z_+(t)\|_{L^\infty}\bigr)+C\bigl(\|W_x(t)\|_{L^\infty}^2+\|W_x(t)\|_{L^\infty}^4\bigr)\nonumber\\
&\quad+CA\bigl(\|R_1(t)\|_{L^2}+\|R_2(t)\|_{L^2}\bigr).
\end{align}
Taking the combination $\eqref{eqs4.1}+\Lambda\eqref{eqs7.1}+M\eqref{eqs5.6}$ and choosing $1\ll \Lambda\ll M$, with $\beta$ sufficiently small in Proposition~\ref{prop7.1} and Lemma~\ref{lem5.1} and $A_0,\delta_0$ sufficiently small, we obtain
\begin{align}\label{eqs8.4}
&E'_{\mathrm{tot}}(t)+c\int_{\R}\bigl[z_-^2(z_-)_x^2+z_+^2(z_+)_x^2+Q^4\bigr]\,dx+c\int_{\R} W_x(z_--z_+)^2\bigl[(z_-+z_+)^2+2z_+^2\bigr]\,dx\nonumber\\
&\quad+c\int_{\R}(z_-^4z_+^2+z_-^2z_+^4)\,dx+cA^2D_1(t)\nonumber\\
&\le k(t)E_{\mathrm{tot}}(t)+CA^4(1+t)^{-3/2}+CE(t)^{3/4}\bigl(\|R_1\|_{L^4}+\|R_2\|_{L^4}+\|[\mu(V)R_1]_x\|_{L^4}\bigr).
\end{align}

By Lemma~\ref{lem5.2} and \eqref{eqs8.2}, $E_{\mathrm{tot}}(t)$ satisfies the cubic Nash inequality
\begin{equation}\label{eqs8.5}
E_{\mathrm{tot}}(t)^3\le CA^8\Bigl(\int_{\R}\bigl[z_-^2(z_-)_x^2+z_+^2(z_+)_x^2+Q^4\bigr]\,dx+A^2D_1(t)\Bigr).
\end{equation}
By \eqref{eqs2.15} and $\min\{1,\delta\sqrt{1+t}\}\le A(1+t)^{1/4}$,
\begin{align}\label{eqs8.6}
\bigl(\|R_1\|_{L^4}+\|R_2\|_{L^4}+\|[\mu(V)R_1]_x\|_{L^4}\bigr)^{4/3}\le CA^{4/3}(1+t)^{-3/2}.
\end{align}
Applying Young's inequality with exponents $(4,4/3)$ gives
\begin{align*}
&E(t)^{3/4}\bigl(\|R_1\|_{L^4}+\|R_2\|_{L^4}+\|[\mu(V)R_1]_x\|_{L^4}\bigr)\\
&\quad\le \beta\frac{E(t)^3}{A^8}+C_\beta A^{8/3}\bigl(\|R_1\|_{L^4}+\|R_2\|_{L^4}+\|[\mu(V)R_1]_x\|_{L^4}\bigr)^{4/3}\\
&\quad\le \beta\frac{E_{\mathrm{tot}}(t)^3}{A^8}+C_\beta A^4(1+t)^{-3/2}.
\end{align*}
Absorbing $\beta A^{-8}E_{\mathrm{tot}}(t)^3$ into the left-hand side of \eqref{eqs8.4} via \eqref{eqs8.5} yields
\begin{align}\label{eqs8.8}
E'_{\mathrm{tot}}(t)+cA^{-8}E_{\mathrm{tot}}(t)^3\le k(t)E_{\mathrm{tot}}(t)+CA^4(1+t)^{-3/2}.
\end{align}

By \eqref{eqs5.4}, \eqref{eqs4.12}, and Lemma~\ref{lem2.1},
\begin{align}\label{eqs8.7}
\int_0^\infty k(t)\,dt\le C+CA^{3/2}+C\sqrt{\delta}+CA\sqrt{\delta}\le C.
\end{align}
Setting $\widetilde E(t):=\exp\bigl(-\int_0^tk(\tau)\,d\tau\bigr)E_{\mathrm{tot}}(t)$, we have $c_1 E_{\mathrm{tot}}(t)\le \widetilde{E}(t)\le E_{\mathrm{tot}}(t)$ for all $t\ge 0$ with $c_1=\exp\bigl(-\int_0^\infty k(\tau)\,d\tau\bigr)>0$, and \eqref{eqs8.8} reduces to
\begin{align}\label{eqs8.9}
\widetilde E'(t)+cA^{-8}\widetilde E(t)^3\le CA^4(1+t)^{-3/2},
\end{align}
with $\widetilde E(0)=E_{\mathrm{tot}}(0)\le CA^4$. Comparing $\widetilde E(t)$ with the barrier $Y(t):=C_1A^4(1+t)^{-1/2}$ for $C_1>0$ chosen large enough that $cC_1^3\ge\frac12C_1+C$ and $C_1A^4\ge\widetilde E(0)$, we obtain
\begin{align*}
E_{\mathrm{tot}}(t)\le CA^4(1+t)^{-1/2}\qquad\forall\,t\ge 0.
\end{align*}

Finally, by \eqref{exact} and $w(v)\ge c_0>0$, we have $|\phi|+|\psi|\le C(|z_-|+|z_+|)$, so $\|(\phi,\psi)(t)\|_{L^4}^4\le CE_{\mathrm{tot}}(t)$, which gives \eqref{eqs2.21} and $\|z_\pm(t)\|_{L^4}\le CA(1+t)^{-1/8}$. Moreover, $A^2E_1(t)\le CE_{\mathrm{tot}}(t)\le CA^4(1+t)^{-1/2}$ together with \eqref{eqs5.5} and \eqref{eqs4.10} implies
\[
\|(\phi_x,h_x)(t)\|_{H^1}\le CA(1+t)^{-1/4} \qquad\text{and}\qquad \|(z_-)_x(t)\|_{L^2}+\|(z_+)_x(t)\|_{L^2}\le CA(1+t)^{-1/4}.
\]
Interpolating $\|z_\pm(t)\|_{L^4}$ with $\|(z_\pm)_x(t)\|_{L^2}$ via the Gagliardo--Nirenberg inequality
\begin{align*}
\|z_-(t)\|_{L^\infty}+\|z_+(t)\|_{L^\infty}\le C\sum_{\pm}\|z_\pm(t)\|_{L^4}^{2/3}\|(z_\pm)_x(t)\|_{L^2}^{1/3}\le CA(1+t)^{-1/6}
\end{align*}
completes the proof of Theorem~\ref{thm2.2}.
\end{proof}

\bigskip
\noindent\textbf{Competing interests.} The authors declare that they have no competing interests.

\medskip
\noindent\textbf{Data availability.} Data sharing is not applicable to this article as no datasets were generated or analyzed during the current study.

\end{document}